\documentclass[a4paper]{article}
\usepackage{amsthm,amsmath,amssymb}
\usepackage{color}

\newtheorem{teor}{Theorem}[section]
\newtheorem{defin}[teor]{Definition}
\newtheorem{lemm}[teor]{Lemma}
\newtheorem{osse}[teor]{Remark}
\newtheorem{prop}[teor]{Proposition}
\newtheorem{defi}[teor]{Definition}
\newtheorem{coro}[teor]{Corollary}
\newtheorem{prob}[teor]{Problem}

\newcommand{\bele}{\begin{lemm}\begin{sl}}
\newcommand{\enle}{\end{sl}\end{lemm}}
\newcommand{\bedef}{\begin{defi}\begin{sl}}
\newcommand{\eddef}{\end{sl}\end{defi}}
\newcommand{\bete}{\begin{teor}\begin{sl}}
\newcommand{\ente}{\end{sl}\end{teor}}
\newcommand{\beos}{\begin{osse}\begin{rm}}
\newcommand{\eddos}{\end{rm}\end{osse}}
\newcommand{\bepr}{\begin{prop}\begin{sl}}
\newcommand{\empr}{\end{sl}\end{prop}}
\newcommand{\bepro}{\begin{prob}\begin{rm}}
\newcommand{\empro}{\end{rm}\end{prob}}
\newcommand{\bede}{\begin{defin}\begin{sl}}
\newcommand{\edde}{\end{sl}\end{defin}}
\newcommand{\beco}{\begin{coro}\begin{sl}}
\newcommand{\enco}{\end{sl}\end{coro}}

\newcommand{\quext}{\quad\text}
\newcommand{\qquext}{\qquad\text}
\newcommand{\de}{\partial}

\newcommand{\RR}{\mathbb{R}}

\newcommand{\beeq}[1]{\begin{equation}\label{#1}}
\newcommand{\eddeq}{\end{equation}}

\newcommand{\beeqa}[1]{\begin{eqnarray}\label{#1}}
\newcommand{\eddeqa}{\end{eqnarray}}

\newcommand{\beal}[1]{\begin{align}\label{#1}}
\newcommand{\eddal}{\end{align}}

\newcommand{\bespl}[1]{\begin{split}\label{#1}}
\newcommand{\edspl}{\end{split}}

\newcommand{\bega}[1]{\begin{gather}\label{#1}}
\newcommand{\edga}{\end{gather}}

\newcommand{\beeqax}{\begin{eqnarray*}}
\newcommand{\eddeqax}{\end{eqnarray*}}

\def\qed{\ifmmode 
  \else \leavevmode\unskip\penalty9999 \hbox{}\nobreak\hfill
  \fi
  \quad\hbox{\hskip.5em\vrule width.4em height.6em depth.05em\hskip.1em}}
\def\endproofsym{\qed}

\def\endnobox{\def\endproofsym{}\end{proof}\def\endproofsym{\qed}}

\newcommand{\no}{\nonumber}

\newcommand{\beeqao}{\begin{eqnarray}\no}
\newcommand{\bealo}{\begin{align}\no}
\newcommand{\besplo}{\begin{split}\no}
\newcommand{\begao}{\begin{gather}\no}

\newcommand{\io}{\int_\Omega}

\newcommand{\iom}{\int_{\Omega_-}}
\newcommand{\iop}{\int_{\Omega_+}}
\newcommand{\iga}{\int_\Gamma}

\newcommand{\epsi}{\varepsilon}
\newcommand{\ee}{_{\varepsilon}}

\newcommand{\OO}{_{\Omega}}
\newcommand{\oo}{_{\Omega}}

\newcommand{\bn}{\boldsymbol{n}}

\newcommand{\dn}{\partial_{\bn}}
\newcommand{\fhi}{\varphi}
\newcommand{\Fhi}{\Phi}

\newcommand{\lhs}{left hand side}
\newcommand{\rhs}{right hand side}

\DeclareMathOperator{\dive}{div}
\DeclareMathOperator{\deriv}{d}

\DeclareMathOperator{\Lip}{Lip}

\newcommand{\HUV}{H^1(0,T;V)}

\newcommand{\LDH}{L^2(0,T;H)}
\newcommand{\LDV}{L^2(0,T;V)}

\newcommand{\LIH}{L^\infty(0,T;H)}
\newcommand{\LIV}{L^\infty(0,T;V)}

\let\TeXchi\chi
\def\chi{{\setbox0 \hbox{\mathsurround0pt
$\TeXchi$}\hbox{\raise\dp0 \copy0 }}}

\newcommand{\calH}{{\mathcal H}}

\newcommand{\calK}{{\mathcal K}}

\newcommand{\calE}{{\mathcal E}}

\newcommand{\dit}{\deriv\!t}
\newcommand{\dis}{\deriv\!s}
\newcommand{\dix}{\deriv\!x}

\newcommand{\dir}{\deriv\!r}

\newcommand{\ddt}{\frac{\deriv\!{}}{\dit}}

\newcommand{\eet}{_{\varepsilon,t}}

\newcommand{\vu}{\boldsymbol{u}}

\newenvironment{bettirev}{\color{blue}}{\color{black}}
\newcommand{\bber}{\begin{bettirev}}
\newcommand{\eber}{\end{bettirev}}

\newenvironment{michelarev}{\color{red}}{\color{black}}
\newcommand{\bmicr}{\begin{michelarev}}
\newcommand{\emicr}{\end{michelarev}}

\newenvironment{giuliorev}{\color{red}}{\color{black}}
\newcommand{\III}{\begin{giuliorev}}
\newcommand{\EEE}{\end{giuliorev}}

\numberwithin{equation}{section}
\begin{document}


\title{On the Cahn-Hilliard-Darcy system with mass source
and strongly separating potential}

\author{Giulio Schimperna\\
Dipartimento di Matematica, Universit\`a di Pavia\\
Via Ferrata~5, I-27100 Pavia, Italy\\
E-mail: {\tt giusch04@unipv.it}
}

\maketitle

\begin{abstract}
 We study an evolutionary system of Cahn-Hilliard-Darcy type
 including mass source and transport effects. The system may arise in a number of physical 
 situations related to phase separation phenomena with convection, with the main 
 and most specific application being related to tumoral processes, where
 the variations of the mass may correspond to growth, or shrinking, of the tumor. 
 We prove existence of weak solutions in the case when the configuration
 potential for the order parameter $\fhi$ is designed in such a way to 
 keep $\fhi$ in between the reference interval $(-1,1)$ despite the 
 occurrence of mass source effects. Moreover, in the two-dimensional
 case, we obtain existence and uniqueness of strong (i.e., more regular)
 solutions.
\end{abstract}


\vspace{2mm}

\noindent {\bf AMS (MOS) subject clas\-si\-fi\-ca\-tion:}%
~~35D30, 35K61, 35Q35, 76D27, 92C30.

\bigskip
\bigskip
\smallskip

\centerline{\sl Dedicated to Maurizio Grasselli on the occasion of his 60th birthday,}

\smallskip

\centerline{\sl with friendship and admiration} 

\medskip


\section{Introduction}

In a bounded domain $\Omega\subset \mathbb{R}^d$, $d=2,3$, with smooth boundary $\Gamma = \partial \Omega$,
we consider the following PDE system:
\begin{align}\label{phasein}
  & \fhi_t + \vu \cdot \nabla \fhi = \Delta \mu + ( 1 - \fhi ) S, \\
 \label{muin}
  & \mu = - \Delta \fhi + f(\fhi) - \lambda \fhi, \\
 \label{Sin}
  & \dive \vu = S, \\
 \label{vuin}
  & \vu = - \nabla p + \mu \nabla \fhi, \\
 \label{S2in}
  & S = - ( 1 + \fhi) \gamma(x,t,\fhi).
\end{align}
Equations \eqref{phasein}-\eqref{muin} are a variant of the well-known 
Cahn-Hilliard system \cite{CH} for phase separation, with $\fhi$ representing the
order parameter, so normalized that the values $\pm1$ correspond to the pure states,
and $\mu$ is an auxiliary variable corresponding to the chemical potential of the 
process. In our setting the phase separation is influenced by transports
effect driven by the movement of the particles of the substance; these effects are
described by means of the macroscopic velocity $\vu$. In particular, the \lhs\ of 
\eqref{phasein} represents the material derivative of $\fhi$. 

The velocity $\vu$ is assumed to satisfy a particular form of the Darcy law, stated by
equations \eqref{Sin} and \eqref{vuin}, where $p$ is the pressure and 
the coupling term $\mu\nabla \fhi$ represents a Korteweg force. It is worth 
noting that other expressions for the pressure and the Korteweg term
could be considered, leading to somehow different models (see the discussion
in Subsec.~\ref{subsec:rema} below). One of the main points in our analysis 
stands in the fact that we assume that there is no conservation of mass; in other
words, the spatial mean value of $\fhi$, corresponding to the proportion between
the two components, or phases, may vary in time. The volumic mass source
depends on the function $S$ on the \rhs\ of \eqref{phasein}, whose expression 
is specified by the last equation \eqref{S2in}. In the sequel we will extensively
comment about this choice; for now we just observe that the
function $\gamma$ is assumed to be measurable, globally bounded, and
Lipschitz continuous with respect to $\fhi$. We admit an explicit dependence
on $(x,t)$ because, in applications, the above system may arise as a subsystem
of a more complex model also involving other quantities (further details
are given below). Finally, the terms $f(\fhi) - \lambda \fhi$ in \eqref{mu} represent
the derivative of what, in the Cahn-Hilliard terminology, is generally
noted as the ``configuration potential'' of the phase separation.
In most cases, such a potential has a double-well structure, 
with the minima attained in proximity of the pure states $\fhi=\pm1$.
In our notation, $F$ is the convex part of the potential, so that $f = F'$, 
while the remainder (concave) part is given by $-\lambda \fhi^2/2$, where
$\lambda\ge0$. Our specific choice for $F$ will be thoroughly discussed in the sequel.

System \eqref{phasein} may have an independent mathematical interest in itself
since it provides a prototypical coupling between the Cahn-Hilliard
equation with mass source and the Darcy law. On the other hand, the choice
of parameters and data (and, in particular, the expression \eqref{S2in}
for the mass source term) is motivated by the specific situation of
diffuse interface models for tumor growth.
Indeed, the recent literature dealing with mathematical models for cancer evolution
is very vast and rapidly expanding (we quote, with no 
claim of completeness, \cite{CL,FGR,FLRS,GL0,GL1,GLSS,JWS,LTZ,OHP}
for an overview).
In tumor models based on the Cahn-Hilliard equation, 
the order parameter $\fhi$ may represent the local concentration
of one class of cells (e.g., ``healthy'', ``cancerous'', or ``dead''), and the presence
of a volumic mass source describes the fact that the tumor may grow, or shrink, depending
on the effects of other factors (like nutrients, drugs, or blood vessels),
here represented by the explicit dependence on $(x,t)$ of the term $S$. 
Here, we will consider a {\sl scalar}\/ order parameter
$\fhi$, which means that we assume that there are only {\sl two}\/ types of cells;
on the other hand, extending our results to the multi-phase case, i.e., to a 
vector-valued variable $\fhi$, should be mathematically straighforward
(see, e.g., \cite{CLLW,FLRS,WLFC} for related work).

Mathematically speaking, the main difficulty of system \eqref{phasein}-\eqref{S2in}
arises because the mass source conflicts with the fact that $\fhi$
represents a normalized order parameter and, for this reason, should take values
into a bounded interval (in our setting $[-1,1]$), whose extrema correspond
to the pure states. More precisely, in tumor applications, $\fhi(x,t) = 1$
may tell us that, at some point $(x,t)$, only healthy cells are present; respectively,
$\fhi(x,t) = - 1$ means that, at $(x,t)$, only cancerous cells occur. Since
the Cahn-Hilliard equation is of the fourth order in space and, 
consequently, its solutions do not obey a maximum principle,
the property $\fhi\in[-1,1]$ is usually enforced by choosing the convex part $F$
of the potential to be ``singular'', i.e.\ by assuming that $F:\RR \to [0,+\infty]$
takes finite values in the interval $[-1,1]$ only,
with the most common example being given by the so-called ``logarithmic
potential'', whose convex part has the expression
\begin{equation}\label{logpot0}
  F(\fhi) = (1 + \fhi)\ln ( 1 + \fhi ) + (1 - \fhi)\ln ( 1 - \fhi ), 
  \qquext{so that }\,\,\,
  f(\fhi) = F'(\fhi) = \ln (1 + \fhi) - \ln (1- \fhi).
\end{equation}
In order to explain how the interaction with the mass source may give troubles, 
we consider the following ``toy situation'', where no transport occurs, $\lambda = 0$,
and the mass source is set equal to the constant value $1$:
\begin{equation}\label{toy}
   \fhi_t - \Delta \mu = 1 , \qquad \mu = - \Delta \fhi + f(\fhi).
\end{equation}
Then, integrating the first equation in space and assuming the (standardly
used for the Cahn-Hilliard equation) no-flux (i.e., homogeneous Neumann)
boundary conditions, one gets $(\fhi\OO)' = 1$ (here and below, $v\OO$ 
denotes the spatial mean of a function, or functional, $v$), so that,
independently of the initial value, at some time $t$ one necessarily
has $\fhi\OO(t) > 1$, which is inconsistent with the fact that, in 
the second equation, $f$ is defined only as its argument ranges 
between $-1$ and $1$. 

To fix this ``inconsistency'' issue, various strategies have been
proposed. Let us now briefly discuss them.
\begin{enumerate}
 \item[{\bf (A)}] Replacing \eqref{logpot0} with a ``smooth'' potential defined over the whole
 real line. For instance one may take $F(\fhi) = \fhi^4$ so that $F(\fhi) - \lambda \fhi^2/2$ 
 is the standard quartic double well potential (cf.~\cite{FGR,GL0,GL1}
 for some examples in the tumor-related literature). This choice, of course, fixes the issue,
 but the interpretation of $\fhi$ as an order parameter is lost since the values of $|\fhi|$ may
 evolve to become arbitrarily large at some points, even when the initial datum
 $\fhi_0$ is everywhere strictly between $-1$ and $1$.
 \item[{\bf (B)}]  Modifying the \rhs\ of \eqref{phasein} in such a way to include an explicit
 dependence on $\mu$ (cf.~\cite{CGH,FGR,Sig} for examples in the tumor-related
 literature). Referring to a ``toy-case'' similar to \eqref{toy}, we obtain 
 something like
 \begin{equation}\label{toy2}
     \fhi_t - \Delta \mu + \mu = 1 , \qquad \mu = - \Delta \fhi + f(\fhi).
 \end{equation}
 Then, replacing the expression of $\mu$ given by the second relation into the 
 first one, one may see that the function $f$ somehow influences the balance of mass, 
 so contributing to keep $\fhi$ in between $-1$ and $1$; on the other hand, the above 
 choice may be questionable from the modeling viewpoint because it
 also adds an ``artificial'' dissipation term in the energy balance law.
 \item[{\bf (C)}]  Designing the mass source term in such a way to constrain the mean value of $\fhi$ in
 between $-1$ and $1$ (cf.~\cite{FLRS,GLRS} for examples in the tumor-related
 literature). Namely, one may take 
 \begin{equation}\label{toy3}
     \fhi_t - \Delta \mu = - K \fhi + s(x,t,\fhi),
 \end{equation}
 where $K>0$ is large enough to dominate the nonlinear part $s$. In this case, integrating 
 \eqref{toy3} in space leads to 
 \begin{equation}\label{toy3m}
    (\fhi\oo)_t + K \fhi\OO = \frac1{|\Omega|} \io s(x,t,\fhi).
 \end{equation}
 At least if the \rhs\ is ``small'' in the $L^\infty$-norm (in a way that may
 depend on $K$ and on the initial datum), using elementary ODE techniques one 
 can easily prove that $\fhi\oo$ is constrained
 between $-1$ and $1$, which is the key step in order to avoid the inconsistency issue mentioned
 above and obtain existence of a weak solution.
\end{enumerate}
We focus here on the approach {\bf (C)}, which was used in the recent
paper \cite{FLRS} dealing with a model very similar to ours.
Actually, when the macroscopic velocity is involved,
a further complication due to the choice of the boundary conditions occurs.
To realize this, we integrate \eqref{phasein} in space:
applying the Gauss-Green formula we then obtain
\begin{equation}\label{medie}
   \ddt \io \fhi = \io S 
    + \iga \dn \mu - \iga \fhi \vu \cdot \bn,
\end{equation}
where $\bn$ is the outer normal unit vector to $\Gamma$. 
Hence, if we take homogeneous Neumann boundary conditions for 
$\mu$, the last integral does {\sl not} disappear and explicitly influences
the balance of mass. Moreover, due to the poor regularity 
expected for $\vu$, it seems to be hard to estimate it.
Hence, obtaining a manageable ODE like \eqref{toy3m}
seems to be out of reach in the no-flux case for $\mu$.
To overcome this issue, in \cite{FLRS}, in place of the Neumann 
condition for $\mu$, was assumed a ``coupled'' boundary condition of the form 
\begin{equation}\label{coupled}
  \dn \mu - \fhi \vu \cdot \bn \equiv 0 \quext{on }\,\Gamma.
\end{equation}
With such a choice, the two boundary integrals in \eqref{medie}
cancel with each other so that a constraint on the spatial mean of $\fhi$ can be
reached provided that the mass source term is suitably designed. 
It is however worth discussing the ``physical'' meaning of condition \eqref{coupled},
which corresponds to prescribing that the mass source effect is 
purely ``volumic'' (due, e.g., to the growth of existing cells and to
the effects of nutrients or drugs), while there is no mass inflow or 
outflow due to transport.

In the present paper, we would like to address the case when the (probably
more natural) no-flux condition for $\mu$ is taken. In such a situation
only the first boundary integral in \eqref{medie} disappears, which means
that there may be some mass inflow or outflow through the boundary due to
transport. For instance, if $\vu$ is oriented as $\bn$ (outflow) and 
$\fhi\sim 1$ (tumor-phase) on $\Gamma$, then the boundary integral
tells us that the (prevailing) cancerous cells are being 
transported away, which is something expected. On the other hand, since the 
boundary integral in general cannot be estimated, in the no-flux case in principle
one cannot exclude that the mean value of $\fhi\OO$ might exceed $1$ (or $-1$), possibly
leading to ``inconsistency'' of the model. Indeed, 
the only certainly available a-priori information comes from
the energy estimate (cf.~\eqref{energy} below), 
which, however, may be not sufficient to guarantee that $\fhi$ stays between $-1$ and $1$. 
In particular, inconsistency may occur in the case of the ``logarithmic'' potential
\eqref{logpot0} because the function $F$ is bounded over the {\sl closed}\/
interval $[-1,1]$ and, for this reason, it does not provide a sufficiently
strong constraint on $\fhi$ in presence of a source of mass combined with
a transport effect. 

In order to avoid this problem, we will consider here what we will term as a
``strongly separating'' potential (cf.~\eqref{logpot} below for the precise
expression), i.e., we will consider a function $F$ which is 
unbounded near $\pm 1$ so to 
provide a stronger ``separation effect'' on $\fhi$ as a mere consequence
of the energy estimate. Note that, from the qualitative point of view, 
the expression of $f=F'$ is rather similar in the two cases~\eqref{logpot}
and~\eqref{logpot0}; hence, this modification of the energy 
does not seem to affect the qualitative aspects of the model. 
On the other hand, we may show that, with \eqref{logpot}, 
global existence of a weak solution can be obtained under very natural assumptions
on the other parameters (and, in particular, on the volumic mass source $S$)
and without adding any ``artificial'' term.

It is worth noting that, even with the help of \eqref{logpot}, the procedure
we use in order to obtain a suitable set of a-priori estimates is nontrivial;
indeed, we still need to explicitly exclude the ``inconsistency'' phenomenon, i.e.,
the fact that $|\fhi\OO|$ may, at some time, exceed $1$. This is
a real issue because the evolution of the mass does not depend, at least directly,
on equation~\eqref{muin} (and, in turn, on the expression of $f$). Nevertheless,
taking advantage of  a proper approximation of the system, we are 
able to rigorously show that an estimate for the spatial mean value 
of $\fhi$ arises just as a direct consequence of the energy bound.

In the proof of existence of weak solutions, we will consider the 
two- and three-dimensional cases together since no variations are necessary
at this level. Of course, one may wonder whether, at least in the two-dimensional 
case, further properties might be proved under additional
assumptions of coefficients and data. This was, indeed, the spirit of 
the result proved in the recent paper \cite{GLRS}. Actually, we will
show that the arguments of \cite{GLRS} can also be adapted to our situation 
so to provide existence and uniqueness
of ``strong'' solutions in 2D; note that this adaptment is not completely
straighforward because the different boundary conditions assumed here
force us to estimate several terms in a different way compared to 
\cite{GLRS}.

\smallskip

The remainder of the paper is organized as follows: in the next section,
we introduce our notation, state our assumptions on coefficients and 
data, and present the precise statement of our results. The proof
of existence mainly relies on a refined version of the energy estimate, which is 
first presented in a formal way in Section~\ref{sec:formal}. Then, 
in the subsequent Section~\ref{sec:appro}, we detail a suitable approximation
of the system and the rigorous version of the a-priori estimates.
Moreover, the limit with respect to the regularization parameter
is achieved by weak compactness argument and a number of additional comments
are given. Finally, in Section~\ref{sec:rego} we present the proof
of the regularity and uniqueness result holding in two space dimensions.


\section{Assumptions and main results}
\label{sec:main}

We assume $\Omega \subset \RR^d$, $d\in\{2,3\}$, be a smooth and bounded domain
with boundary $\Gamma=\de\Omega$. In the existence proof, 
embeddings, interpolation inequalities, and related exponents,
will be used by referring to the case $d=3$; of course in the 
two-dimensional case the results remain valid and, in addition, the regularity 
properties of solutions may be improved.

We set $H=L^2(\Omega)$, $V:=H^1(\Omega)$ and $V_0:=H^1_0(\Omega)$. For notational
simplicity we will use the same letters $H$, $V$ and $V_0$ to denote vector-valued
spaces (for instance, $H$ may also refer to $L^2(\Omega)^d$). Then, identifying
$H$ with $H'$ by means of the standard $L^2$-scalar product, we obtain the Hilbert
triplets $V\subset H\subset V'$ and $V_0\subset H\subset V_0'$. We will denote by
$\| \cdot \|$ the norm in $H$ and by $\| \cdot \|_X$ the norm in a generic Banach
space $X$. We also indicate by $A : V\to V'$ and $B:V_0 \to V_0'$ the weak version of 
the Neumann and of the Dirichlet Laplacian, respectively.

We consider the system
\begin{align}\label{phase}
  & \fhi_t = \Delta \mu + ( 1 - \fhi ) S - \vu \cdot \nabla \fhi, \\
 \label{mu}
  & \mu = - \Delta \fhi + f(\fhi) - \lambda \fhi, \\
 \label{S}
  & \dive \vu = S, \\
 \label{vu}
  & \vu = - \nabla p + \mu \nabla \fhi, \\
 \label{S2}
  & S = - ( 1 + \fhi) \gamma(x,t,\fhi),
\end{align}
coupled with no-flux conditions for $\fhi$ and $\mu$ and with homogeneous Dirichlet
conditions for $p$. Namely, we assume
\begin{equation}\label{bound}
  \dn \mu = \dn \fhi = p = 0 \quext{on }\,\Gamma.
\end{equation}
As said, $\lambda \ge 0$ is a given constant and we will assume that the 
convex part $F$ of the configuration potential takes the following
``strongly separating'' expression:
\begin{equation}\label{logpot}
  F(\fhi) = - \ln ( 1 - \fhi^2 ), \qquext{so that }\,\,\,
  f(\fhi) = F'(\fhi) = \frac{2\fhi}{1-\fhi^2}.
\end{equation}
Note that the above is just a prototypical choice; more general forms
of $F$ could indeed be taken. The key point stands in the fact that, differently
from the case \eqref{logpot0} of the ``standard'' Cahn-Hilliard logarithmic
potential, here not only the derivative $f$ is singular at $\pm1$, but
the same is true also for $F$. As noted in the introduction, this property
is crucial in order to avoid the ``inconsistency'' phenomenon due to 
the mass source effect.

The function $\gamma$ characterizing the mass source term is 
assumed to satisfy 
\begin{align}\label{hp:gamma}
  & \gamma\in L^\infty(\Omega\times(0,T)\times\RR),
   \qquad \gamma(x,t,\cdot) \in \Lip(\RR)~~\text{for a.e.~}\,(x,t)\in \Omega\times(0,T),\\
 \label{hp:gamma2}
  & \gamma(x,t,r) \equiv 0~~\text{for all }\,|r|\ge 2~~
   \text{and a.e.~}\,(x,t)\in \Omega\times(0,T).
\end{align}
It is worth commenting a bit the above assumptions. First of all, since $\gamma$ represents
a source of mass, assuming global boundedness of it is a natural condition; moreover
hypothesis \eqref{hp:gamma2} serves just as a normalization property for the sake
of building a sound approximation of the system; indeed, 
in the limit, $\fhi$ will take its values in the interval $(-1,1)$;
hence the behavior of $\gamma$ for $\fhi$ outside $(-1,1)$ is 
factually irrelevant. The Lipschitz continuity with respect to $\fhi$ stated 
by the second \eqref{hp:gamma} is a natural requirement as we aim to 
apply a local existence result to a Faedo-Galerkin regularization.

Looking at the \rhs\ of \eqref{phase}, it is worth observing that
\begin{equation}\label{gamma:2}
  (1-\fhi)S  = - ( 1 - \fhi^2 ) \gamma (x,t,\fhi);
\end{equation}
In a sense this prescribes that, as far as $\fhi$ is close to $1$, or to $-1$, 
the volumic mass source tends to approach $0$ (recall that $\gamma$ is bounded).
This ansatz is rather standard at least in tumor applications,
meaning for instance that, when the tumor cells are strongly
prevailing (compared to the healthy ones), very few new ones 
can be created (see also Remark~\ref{rem:gamma} below for further comments).

Moreover, we complement the system with the initial condition
\begin{equation}\label{init}
  \fhi|_{t=0} = \fhi_0,
\end{equation}
where the initial datum $\fhi_0$ is assumed to satisfy the following regularity
properties, which basically correspond to the finiteness of the physical energy at the initial
time:
\begin{equation}\label{hp:init}
  \fhi_0 \in V, \qquad F(\fhi_0) \in L^1(\Omega),
\end{equation}
Note that then, by Jensen's inequality, there follows
\begin{equation}\label{jen}
  F \bigg( \io \fhi_0(x)\,\frac{\dix}{|\Omega|} \bigg) 
   \le \io F(\fhi_0(x))\,\frac{\dix}{|\Omega|} 
   < + \infty;
\end{equation}
Consequently, due to the choice \eqref{logpot}, and more specifically 
to the fact that $\lim_{|r|\to 1} F(r) = +\infty$, one has 
\begin{equation}\label{hp:init2}
  (\fhi_0)\OO \in (-1,1).
\end{equation}
This fact is crucial in order to get a control of the spatial mean of $\fhi$ over the 
interval $(0,T)$. We notice that, in the case of the logarithmic potential \eqref{logpot0}, 
\eqref{hp:init2} generally needs to be taken as an additional assumption
(since $F$ is finite over the {\sl closed}\/ interval $[-1,1]$ in that case),
while in our case is just a direct consequence of the expression \eqref{logpot}.

We can now state the main result of this paper:
\bete\label{teo:main}
 Let assumptions\/ \eqref{logpot},
 \eqref{hp:gamma}-\eqref{hp:gamma2}, and\/ \eqref{hp:init} hold. Then,
 there exists at least one quadruplet $(\fhi,\mu,p,\vu)$ of functions such that
 \begin{align} \label{rego:fhi}
  & \fhi \in H^1(0,T;V'+L^1(\Omega)) \cap \LIV \cap L^4(0,T;H^2(\Omega)) \cap L^2(0,T;W^{2,6}(\Omega)),\\ 
  \label{rego:ffhi}
  & f(\fhi) \in L^2(0,T;L^6(\Omega)),\\
  \label{rego:mu}
  & \mu \in \LDV,\\
  \label{rego:vu}
  & \vu \in \LDH,\\ 
  \label{rego:p}
  & p \in L^2(0,T;W^{1,3/2}_0(\Omega)),
\end{align}
 satisfying system\/ \eqref{phase}-\eqref{S2} with the boundary conditions \eqref{bound}
 in the following weak form:
 \begin{align}\label{phaseteo}
   & \fhi_t + A \mu = ( 1 - \fhi ) S - \vu \cdot \nabla \fhi \quext{in } L^2(0,T;V' + L^1(\Omega)), \\
  \label{muteo}
   & \mu = - \Delta \fhi + f(\fhi) - \lambda \fhi, \quext{a.e.~in }\,\Omega\times (0,T), \\
  \label{bounfhi}
   & \dn \fhi = 0, \quext{a.e.~on }\,\Gamma\times (0,T), \\
  \label{Steo}
   & \io \vu(t) \cdot \nabla \eta = - \io S(\cdot,t,\fhi) \eta  \quext{for every }\,\eta \in V_0~~\text{and a.e.\ }\, t \in (0,T), \\
  \label{vuteo}
   & \vu = - \nabla p + \mu \nabla \fhi, \quext{a.e.~in }\,\Omega\times (0,T), \\
  \label{S2teo}
   & S = - ( 1 + \fhi) \gamma(x,t,\fhi), \quext{a.e.~in }\,\Omega\times (0,T),
 \end{align}
 and complying with the initial condition\/ \eqref{init}.
\ente
\noindent%
We conclude this section by detailing our second result, which is devoted to existence and
uniqueness of ``strong'' (i.e.\ more regular) solutions in two dimensions of 
space, so extending to the present setting \cite[Theorem~4.1]{GLRS}.
\bete\label{teo:rego}
 Let assumptions\/ \eqref{logpot},
 \eqref{hp:gamma}-\eqref{hp:gamma2}, and\/ \eqref{hp:init} hold and
 let $\Omega\subset \RR^2$. Moreover, let us set 
 \begin{equation}\label{defi:mu0}
   \mu_0 := - \Delta \fhi_0 + f(\fhi_0) - \lambda \fhi_0
 \end{equation}
 and, correspondingly, let us assume 
 \begin{align}\label{init:reg1}
   & \fhi_0 \in H^2(\Omega), \qquad f(\fhi_0) \in H, \qquad 
    \dn \fhi_0 = 0~~\text{on }\, \Gamma,\\
  \label{init:reg2}
   & \mu_0 \in V.
 \end{align}
 Then, there exists one {\rm and only one} solution $(\fhi,\mu,p,\vu)$ to 
 system\/ \eqref{phase}-\eqref{S2} with the boundary conditions \eqref{bound}
 in the following regularity class:
 \begin{align} \label{rego:fhis}
  & \fhi \in W^{1,\infty}(0,T;V') \cap \HUV \cap L^\infty(0,T;W^{2,r}(\Omega))
   \quext{for all }\,r\in[1,\infty),\\ 
  \label{rego:ffhis}
  & f(\fhi) \in L^\infty(0,T;L^r(\Omega))
   \quext{for all }\,r\in[1,\infty),\\ 
  \label{rego:mus}
  & \mu \in \LIV \cap L^4(0,T;H^2(\Omega)),\\
  \label{rego:vus}
  & \vu \in L^\infty(0,T;H^1(\Omega)),\\ 
  \label{rego:ps}
  & p \in L^\infty(0,T;H^2(\Omega) \cap H^1_0(\Omega)).
 \end{align}
\ente
\noindent%
Note that the exponent $6$ occurring in \eqref{rego:fhi}-\eqref{rego:ffhi} has been replaced
by $r\in[1,\infty)$ in \eqref{rego:fhis}-\eqref{rego:ffhis}. This is 
a consequence of better 2D embeddings (and the same could be done in 
Theorem~\ref{teo:main} when restricted to the two-dimensional setting).

The proof of the above results will occupy the remainder of the paper.
In particular, in the proof of Theorem~\ref{teo:rego}, presented in Section~\ref{sec:rego}
below, we will mainly focus on the regularity part of the statement because uniqueness 
works very similarly with \cite{GLRS} and, for this reason, will be only
sketched.


\section{Formal energy estimate}
\label{sec:formal}

In order to fix the main points of our procedure and understand the role of
the assumptions, we derive, for the reader's convenience,
a formal version of the energy estimate. This means that we will directly
work on the system \eqref{phase}-\eqref{S} without referring to any
approximation. The procedure will be made rigorous in
the next section, where a regularization of the system will be proposed.
We decided to present both versions of the estimate because the rigorous
estimate a bit more technical.

That said, we test \eqref{phase} by $\mu$ and \eqref{mu} by $\fhi_t$ to get
\begin{equation}\label{energy}
  \ddt \calE 
   + \| \nabla \mu \|^2
   = \io (1 - \fhi) S \mu - \io \vu \cdot \nabla \fhi \mu,
\end{equation}
where $\calE$ denotes the standard Cahn-Hilliard energy, i.e., 
\begin{equation}\label{defiE}
  \calE = \calE(\fhi) = \frac12 \| \nabla \fhi \|^2
   + \io \Big( F(\fhi) - \frac\lambda2 \fhi^2 \Big).
\end{equation}
Now, using the expressions \eqref{mu} for $\mu$ and
\eqref{S2} for $S$, it is easy to see that
\begin{equation}\label{co10}
  \io (1 - \fhi) S \mu 
   = - \io (1 - \fhi^2) \gamma(x,t,\fhi) \big( - \Delta\fhi + f(\fhi) - \lambda\fhi \big).
\end{equation}
Then, using \eqref{hp:gamma2} and the fact that $\fhi$ takes values in $(-1,1)$, 
we have
\begin{equation}\label{co11}
  - \io (1 - \fhi^2) \gamma(x,t,\fhi) \big( - \Delta\fhi - \lambda\fhi \big)
   \le c + \frac14 \| \Delta \fhi \|^2.
\end{equation}
On the other hand, using the second \eqref{logpot} and the fact
$\fhi\in(-1,1)$, we deduce
\begin{equation}\label{co11b}
  - \io (1 - \fhi^2) \gamma(x,t,\fhi) f(\fhi)
  = - 2 \io \gamma(x,t,\fhi) \fhi
  \le c.
\end{equation}
\beos\label{rem:gamma}
 In \eqref{co11b} we used the degenerate behavior of $(1-\fhi)^2$ at $\pm 1$ in order to
 compensate the singular character of $f(\fhi)$. Note that the same argument would
 work also in the case of the standard ``logarithmic'' potential \eqref{logpot0}.
 On the other hand, for different expressions of the mass source term (for instance
 if one forgets the factor $(1+\fhi)$ in the expression \eqref{S2} for $S$), the 
 argument may fail. However, an estimate could still be obtained if
 additional sign conditions on $\gamma$ are assumed.
 For instance, assuming for simplicity $\gamma$ to depend 
 only on $\fhi$, if $\gamma(\fhi)$ has the same sign as $\fhi$
 at least for $|\fhi|$ close to $1$,
 then the integral in \eqref{co11b} may be moved to the \lhs\
 and gives a positive contribution (so it does not need to be controlled). On the 
 other hand, we preferred to avoid sign conditions on $\gamma$ because in 
 applications its expression (and in particular its sign) may be determined by
 the effects of other quantities (cf., for instance, the discussion in 
 \cite[Sec.~3.4]{MRS}). We may also observe that the expression \eqref{S2}
 for $S$ appears to be realistic at least in tumor-related applications 
 (see, e.g., \cite{GLSS}).
\eddos
\noindent%
Next, testing \eqref{vu} by $\vu$, integrating by
parts, and using \eqref{S} with the boundary conditions
\eqref{bound}, we obtain
\begin{equation}\label{co12}
  \| \vu \|^2 
  = \io p S + \io \vu \cdot \nabla \fhi \mu.
\end{equation}
Hence, summing \eqref{co12} to \eqref{energy} and using \eqref{co11}-\eqref{co11b}
we infer
\begin{equation}\label{co13}
  \ddt \calE 
   + \| \nabla \mu \|^2
   + \| \vu \|^2 
  \le \io p S 
  + c 
  + \frac14 \| \Delta \fhi \|^2.
\end{equation}
To estimate the integral term on the \rhs, we 
define $\zeta$ as the solution of the 
following time-dependent family of elliptic problems:
\begin{equation}\label{co14}
  B \zeta = S. 
\end{equation}
Taking the divergence of \eqref{vu}, and recalling
the boundary conditions, we can write
\begin{equation}\label{co15}
  B p = S - \dive( \mu \nabla \fhi ).
\end{equation}
Testing \eqref{co15} by $\zeta=B^{-1} S$ and using the boundary conditions, 
we infer
\begin{align}\no
  \io p S 
   & = \| S \|_{H^{-1}(\Omega)}^2 - \io \dive( \mu \nabla \fhi )\zeta \\
 \no
   & = \| S \|_{H^{-1}(\Omega)}^2 + \io \mu \nabla \fhi \cdot \nabla\zeta \\
 \no
   & = \| S \|_{H^{-1}(\Omega)}^2 
    + \io \big( - \Delta \fhi + f(\fhi) - \lambda \fhi \big) \nabla \fhi \cdot \nabla \zeta \\
 \no
   & = \| S \|_{H^{-1}(\Omega)}^2 
    - \io \Delta \fhi \nabla \fhi \cdot \nabla \zeta
    + \io f(\fhi) \nabla \fhi \cdot \nabla \zeta
    - \io \lambda\fhi \nabla \fhi \cdot \nabla \zeta\\
 \label{co16}
  & =: \| S \|_{H^{-1}(\Omega)}^2 
    + I_1 + I_2 + I_3
    \le c + I_1 + I_2 + I_3,
\end{align}
the last inequality following from the fact that $S$ is bounded in the 
$L^\infty$-norm.

We now provide a control of the integral terms on the \rhs. Firstly, 
by elementary use of interpolation and Sobolev's embeddings,
we have
\begin{align}\no
  I_1 
   & = - \io \Delta \fhi \nabla \fhi \cdot \nabla \zeta
   \le \| \Delta \fhi \| \| \nabla \fhi \|_{L^3(\Omega)} \| \nabla \zeta \|_{L^6(\Omega)}\\
 \label{co17}
   & \le \| \Delta \fhi \|^{3/2} \| \nabla \fhi \|^{1/2} \| \Delta \zeta \|
    \le \frac14 \| \Delta \fhi \|^2 + c \| \nabla \fhi \|^2,
\end{align}
where, to deduce the last inequality, we also used \eqref{co14} and 
the uniform boundedness of $S$.

Second, we observe that
\begin{align}\no
  I_2 
   & = \io f(\fhi) \nabla \fhi \cdot \nabla \zeta
    \le \| f(\fhi) \nabla \fhi\|_{L^1(\Omega)} \| \nabla \zeta\|_{L^\infty(\Omega)} \\
 \no
   &  \le \| f(\fhi) \nabla \fhi\|_{L^1(\Omega)} \| B \zeta\|_{L^4(\Omega)} 
    \le \| f(\fhi) \nabla \fhi\|_{L^1(\Omega)} \| S \|_{L^4(\Omega)} \\
 \label{co18}
   &  \le c \| f(\fhi) \nabla \fhi\|_{L^1(\Omega)} \| S \|_{L^\infty(\Omega)}
   \le c \| f(\fhi) \nabla \fhi\|_{L^1(\Omega)}.
\end{align}
Finally, using again that $\fhi\in(-1,1)$, we have
\begin{equation}\label{co18b}
  I_3 = - \io \lambda\fhi \nabla \fhi \cdot \nabla \zeta
  \le c \| \nabla \fhi \|
  \le c \big( 1 + \| \nabla \fhi \|^2 \big).
\end{equation}
To control the \rhs s of \eqref{co17} and \eqref{co18}, the {\sl strongly separating}\/ behavior
of $F$ is crucial. To exploit it, we actually need to perform a further calculation by 
testing \eqref{mu} by $-\Delta \fhi$ to obtain
\begin{equation}\label{co19}
  \| \Delta \fhi \|^2 
   + \io f'(\fhi) |\nabla \fhi|^2 
    \le \frac12 \| \nabla \mu \|^2 +  \Big( \frac12 + \lambda \Big) \| \nabla \fhi \|^2.
\end{equation}
Summing the above to \eqref{co13} and using \eqref{co17} and \eqref{co18},
we then infer
\begin{align}\no
  & \ddt \calE 
   + \frac12 \| \nabla \mu \|^2
   + \| \vu \|^2 
   + \frac12 \| \Delta \fhi \|^2 
   + \io f'(\fhi) |\nabla \fhi|^2 \\
 \label{co1a}
  & \mbox{}~~~~~   
   \le  c \| f(\fhi) \nabla \fhi\|_{L^1(\Omega)}
   + c \big( 1 + \| \fhi \|_V^2 \big).
\end{align}
Now, from \eqref{logpot} we may compute
\begin{equation}\label{logpot2}
  f'(\fhi) = \frac{2(1+\fhi^2)}{(1-\fhi^2)^2}.
\end{equation}
Hence, it is easy to see that 
\begin{equation}\label{co1b}
  c \| f(\fhi) \nabla \fhi\|_{L^1(\Omega)}
   = \io \frac{2|\fhi|}{1-\fhi^2} | \nabla \fhi |
   \le \frac12 \io f'(\fhi) |\nabla \fhi|^2 + c.
\end{equation}
Substituting \eqref{co1b} into \eqref{co1a}, we finally arrive at the 
relation
\begin{equation}\label{co1d}
 \ddt \calE 
   + \frac12 \| \nabla \mu \|^2
   + \| \vu \|^2 
   + \frac12 \| \Delta \fhi \|^2 
   + \frac12 \io f'(\fhi) |\nabla \fhi|^2 
   \le c \big( 1 + \| \fhi \|_V^2 \big).
\end{equation}
Using the Gr\"onwall lemma and noting that, by assumptions~\eqref{hp:init},
the energy $\calE$ is finite at the initial time, it is possible to deduce 
a global estimate. In the next section we will see that, in fact, such
a procedure may be adapted to work in the framework of a rigorous approximation.


\section{Approximation}
\label{sec:appro}

In this part we introduce an approximation of system \eqref{phase}-\eqref{S2}
and sketch a way to prove existence to the obtained regularized system.

Given $\epsi\in (0,1/4)$ which will be let go to $0$
in the limit, the main point stands in providing a suitable regularization
$f\ee$ of the function $f$ compatible with the a-priori estimates. To this aim, 
we take $f\ee$ equal to $f$ in the interval $[-1+\epsi,1-\epsi]$ and $f\ee$ given by the first 
order Taylor expansion of $f$ centered in $1-\epsi$ in $(1-\epsi,+\infty)$
and, respectively, by the first order Taylor expansion of $f$ centered in $-1+\epsi$ in 
$(-\infty,-1+\epsi)$. Then, noting that
\begin{align}\label{fee:11}
  & F(1-\epsi) = F(-1+\epsi) = - \ln (2 - \epsi) - \ln (\epsi) \sim - \ln \epsi, \\
 \label{fee:12}
  & f(1-\epsi) = \frac{2 - 2\epsi}{2 \epsi - \epsi^2} \sim \frac1\epsi,
  \qquad f(- 1 + \epsi) = \frac{-2 + 2\epsi}{2 \epsi - \epsi^2} \sim - \frac1\epsi,\\
 \label{fee:13}
  & f'(1-\epsi) = f'(-1+\epsi) = \frac{2 ( 1 + (1 - \epsi)^2 )}{( 1- (1-\epsi)^2)^2} \sim \frac1{\epsi^2},
\end{align}
and restricting our attention to the case $\fhi>1-\epsi$ (the case $\fhi<-1+\epsi$
being analogous), we have 
\begin{align}\label{fee:21}
  & f\ee(\fhi) = f(1-\epsi) + f'(1-\epsi) (\fhi - 1 + \epsi) 
   \sim \frac1\epsi + \frac1{\epsi^2}(\fhi - 1 + \epsi),\\
 \label{fee:22}
  & f\ee'(\fhi) = f'(1-\epsi) \sim \frac1{\epsi^2},
\end{align}
so that, still for $\fhi>1-\epsi$,  we also obtain  
\begin{equation}\label{fee:23}
  F\ee(\fhi) = F(1-\epsi) + \int_{1-\epsi}^\fhi f\ee(r)\,\dir
   \sim - \ln \epsi +  \frac1\epsi (\fhi - 1+ \epsi) + \frac1{2\epsi^2}(\fhi - 1+ \epsi)^2.
\end{equation}
We then replace $f$ with $f\ee$ in equation \eqref{mu}. Since $f\ee$ is globally Lipschitz, 
at the approximate level the order parameter will no longer
take its values into $(-1,1)$. For this reason, we also need to truncate the terms
$(1-\fhi)$ in \eqref{phase} and $(1+\fhi)$ in \eqref{S2} by replacing them with their
positive parts, noted by $(\cdot)^+$ in the sequel.
Then, in order to prepare the ground for an existence result, we 
rewrite the system by eliminating $\vu$ and using instead the pressure which is a somehow
more natural variable because it satisfies an elliptic equation with proper
(namely, Dirichlet) boundary conditions.

By the above considerations, noting as $(\fhi\ee,\mu\ee,p\ee)$ a candidate
solution, we are led to the system
\begin{align}\label{phaseee}
  & \fhi\eet + A \mu\ee = ( 1 - \fhi\ee )^+ S\ee + \nabla p\ee \cdot \nabla \fhi\ee
   - \mu\ee | \nabla \fhi\ee |^2,\\
 \label{muee}
  & \mu\ee = A \fhi\ee + f\ee(\fhi\ee) - \lambda\fhi\ee, \\
 \label{pee}
  & B p\ee = S\ee - \dive (\mu\ee \nabla\fhi\ee),\\
 \label{S2ee}
  & S\ee = - ( 1 + \fhi\ee)^+ \gamma(x,t,\fhi\ee),
\end{align}
where the boundary conditions are incorporated into the operators $A$ and $B$.

The above system is naturally complemented with the initial conditions \eqref{init},
where the initial data do not need to be regularized. Then, in order to prove 
a local existence result via a discretization scheme, it
is also convenient to eliminate the auxiliary variable $\mu\ee$. Actually,
equations \eqref{phaseee} and \eqref{muee} can be combined as 
\begin{equation}\label{phasemuee}
 \fhi\eet + A^2 \fhi\ee + A f\ee(\fhi\ee) - \lambda A \fhi\ee 
   = ( 1 - \fhi\ee )^+ S\ee + \nabla p\ee \cdot \nabla \fhi\ee
   - (A \fhi\ee + f\ee(\fhi\ee) - \lambda\fhi\ee)| \nabla \fhi\ee |^2.
\end{equation}
Analogously, one may insert \eqref{muee} into \eqref{pee} in order
to eliminate $\mu\ee$. This yields
\begin{equation}\label{pmuee} 
  B p\ee = S\ee - \dive \big( (A \fhi\ee + f\ee(\fhi\ee) - \lambda\fhi\ee)\nabla\fhi\ee\big).
\end{equation}
The resulting system \eqref{phasemuee}-\eqref{pmuee}, with $S\ee$ specified by 
\eqref{S2ee}, can then be solved, at least locally in time, by using the Faedo-Galerkin scheme
and possibly implementing a fixed point argument. This procedure works similarly with
other Cahn-Hilliard-based systems and is therefore omitted. Note that one may
need to use two families of eigenfunctions as Faedo-Galerkin bases because
$\fhi\ee$ and $p\ee$ satisfy different types of boundary conditions.

This procedure yields a local in time solution $(\fhi\ee,p\ee)$. Of course, once such a solution
is obtained, one can go back to the formulation \eqref{phaseee}-\eqref{S2ee} by
defining $\mu\ee$ in the natural way. In addition to that, we note that,
in view of the fact that the a-priori estimates derived below have a global in time 
character, by standard extension arguments one may prove that, in fact, the limit 
solution will be globally defined. For the sake of simplicity, 
we shall directly work on $(0,T)$ also at the approximated level leaving the details 
of the extension argument to the reader.


\subsection{Rigorous a priori estimates}
\label{subsec:compa}

Here, we would like to show that, once one tries to adapt the a-priori estimates
of Section~\ref{sec:formal}
to the regularized system (written in the ``extended'' form \eqref{phaseee}-\eqref{S2ee}),
the procedure still remains valid. To see this, we start with repeating the energy estimate \eqref{energy} in
the present setting. This corresponds, of course, to testing \eqref{phaseee} by
$\mu\ee$, \eqref{muee} by $\fhi\eet$, and \eqref{pee} by $p\ee$, so to obtain
\begin{align}\no
  & \ddt \calE\ee
   + \| \nabla \mu\ee \|^2
   + \| \nabla p\ee \|^2 + \io \mu\ee^2 |\nabla \fhi\ee|^2 
    - 2 \io \mu\ee \nabla\fhi\ee \cdot \nabla p\ee\\
 \label{energyee}
  & \mbox{}~~~~~
   = \io (1 - \fhi\ee)^+ S\ee \mu\ee 
   + \io S\ee p\ee,
\end{align}
where $\calE\ee$ denotes the regularized energy, i.e.,
\begin{equation}\label{defiEee}
  \calE\ee = \calE\ee( \fhi\ee) = \frac12 \| \nabla \fhi\ee \|^2
   + \io \Big( F\ee(\fhi\ee) - \frac\lambda2 \fhi\ee^2 \Big).
\end{equation}
Setting back
\begin{equation}\label{backuee}
  \vu\ee:= - \nabla p\ee + \mu\ee \nabla \fhi\ee,
\end{equation}
we then observe that 
\begin{equation}\label{uee}
  \| \vu\ee \|^2 = \| \nabla p\ee \|^2 + \io \mu\ee^2 |\nabla \fhi\ee|^2 
    - 2 \io \mu\ee \nabla\fhi\ee \cdot \nabla p\ee.
\end{equation}
As before, we can also add to \eqref{energyee}
the contribution of \eqref{muee} tested by $A\fhi\ee$,
i.e., the analogue of \eqref{co19}. Performing standard manipulations, 
it is then not difficult to arrive at
\begin{align}\no
  & \ddt \calE\ee
   + \frac12 \| \nabla \mu\ee \|^2
   + \| A \fhi \ee \|^2 
   + \| \vu \ee \|^2 
   + \io f\ee'(\fhi\ee) | \nabla \fhi\ee |^2\\
 \label{energyee1}
  & \mbox{}~~~~~
   \le \io (1 - \fhi\ee)^+ S\ee \mu\ee 
   + \io S\ee p\ee
   + c \| \nabla \fhi\ee \|^2,
\end{align}
and we need to manage the \rhs. First of all, we can treat the integral of $S\ee p\ee$
similarly as before. Namely, we obtain the analogue of the contributions
$I_j$, $j=1,2,3$, of \eqref{co16}, which we need to control.
The estimate \eqref{co17} of $I_1$ can be repeated
without any variation. Concerning $I_3$, since we do not know the 
a-priori boundedness of $\fhi\ee$ at this level,
we can modify \eqref{co18b} as follows:
\begin{equation}\label{co18bee}
  I_3 = - \io \lambda\fhi\ee \nabla \fhi\ee \cdot \nabla \zeta\ee
  \le c \| \fhi\ee \|_{L^3(\Omega)} \| \nabla \fhi\ee \| \| \nabla \zeta\ee \|_{L^6(\Omega)} 
  \le c \| \fhi\ee \|_{V}^2,
\end{equation}
where standard embeddings have been used.

Actually, the main difference regards the control of 
the term $I_2$ in \eqref{co16}.
Repeating \eqref{co18} with notational variations only, we 
get the $L^1$-norm of $f\ee(\fhi\ee)\nabla\fhi\ee$ on the \rhs, and we
would like to adapt \eqref{co1b} in order to estimate it.
To this aim, for a.e.~$t\in (0,T)$,
we may split $\Omega = \Omega_- \cup \Omega_0 \cup \Omega_+$, where 
\begin{equation}\label{splitomega} 
  \Omega_0 = \Omega_0(t) = \big\{x\in \Omega:~|\fhi\ee(x,t)|\le 1-\epsi\big\}, \qquad 
   \Omega_+ = \Omega_+(t) = \big\{x\in \Omega:~ \fhi\ee(x,t) > 1-\epsi\big\},
\end{equation}
and $\Omega_-$ defined similarly (the dependence on $\epsi$ of the subdomains
is not stressed in the notation). Then, the integral on the \lhs\ of \eqref{co1b}
is decomposed into its components on the three subdomains. 
Clearly, the part on $\Omega_0$ may be treated
as in the previous section, while we need to focus on the components on $\Omega_-$
and $\Omega_+$, and, for brevity, we will only consider the latter.
Actually, using \eqref{fee:21} and (twice) Young's inequality, we obtain
\begin{align}\no
  \iop f\ee(\fhi\ee) | \nabla \fhi\ee |
   & \sim \frac1\epsi \iop  | \nabla \fhi\ee |
    + \frac1{\epsi^2} \iop (\fhi\ee - 1 + \epsi) | \nabla \fhi\ee |\\
 \label{ee:12}
   & \le c + \frac1{4\epsi^2} \iop  | \nabla \fhi\ee |^2
    + \frac1{2\epsi^2} \iop (\fhi\ee - 1 + \epsi)^2
    + \frac1{2\epsi^2} \iop | \nabla \fhi\ee |^2,
\end{align}
where we notice that the constant $c$ may be taken independent of $\epsi$. 

In order to control the above \rhs\ uniformly with respect to $\epsi$, we 
take advantage of the last term on the \lhs\ of \eqref{energyee1}, which, 
owing to \eqref{fee:22}, can be rewritten as follows:
\begin{equation}\label{ee:31}
   \io f\ee'(\fhi\ee) |\nabla \fhi\ee|^2 
    \sim \int_{\Omega_0} f'(\fhi\ee) |\nabla \fhi\ee|^2 
    + \frac1{\epsi^2} \iop |\nabla \fhi\ee|^2 
    + \frac1{\epsi^2} \iom |\nabla \fhi\ee|^2. 
\end{equation}
We also notice that there exists a constant $k>0$, depending on $\lambda$ but
independent of $\epsi$, such that
\begin{equation}\label{coercEee}
  \Phi\ee(r) := \frac12 F\ee(r) - \frac\lambda2 r^2 + k \ge 0
\end{equation}
for every $\epsi\in(0,1/4)$ and $r\in \RR$. 
Using the above relations, \eqref{energyee1} gives
\begin{align}\no
  & \ddt \bigg[ \frac12 \| \nabla \fhi\ee \|^2
   + \io \Phi\ee(\fhi\ee)  
   + \frac12 \io F\ee(\fhi\ee) \bigg]
   + \frac12 \| \nabla \mu\ee \|^2
   + \frac34 \| A \fhi\ee \|^2
   + \| \vu\ee \|^2\\
  \no
 & \mbox{}~~~~~
   + \frac12 \int_{\Omega_0} f'(\fhi\ee) |\nabla \fhi\ee|^2 
   + \frac{1}{4\epsi^2} \iop | \nabla \fhi\ee|^2
   + \frac{1}{4\epsi^2} \iom | \nabla \fhi\ee|^2
    \lesssim \io (1 - \fhi\ee)^+ S\ee \mu\ee 
   + c\\
 \label{energyee3}
  & \mbox{}~~~~~~~~~~
   + \frac{1}{2\epsi^2} \iop (\fhi\ee - 1 + \epsi)^2
   + \frac{1}{2\epsi^2} \iom (\fhi\ee + 1 - \epsi)^2
   + c \| \fhi\ee \|_V^2.
\end{align}
Now, let us observe that, thanks to \eqref{fee:23}, there holds
\begin{equation}\label{ee:32}
  \frac{1}{2\epsi^2} \iop (\fhi\ee - 1 + \epsi)^2
  + \frac{1}{2\epsi^2} \iom (\fhi\ee + 1 - \epsi)^2
   \lesssim \io F\ee(\fhi\ee).
\end{equation}
Hence, \eqref{energyee3} can be rewritten in the simpler form
\begin{align}\no
  & \ddt \bigg[ \frac12 \| \nabla \fhi\ee \|^2
   + \io \Phi\ee(\fhi\ee)  
   + \frac12 \io F\ee(\fhi\ee) \bigg]
   + \frac12 \| \nabla \mu\ee \|^2
   + \frac34 \| A \fhi\ee \|^2
   + \| \vu\ee \|^2\\
  \no
 & \mbox{}~~~~~~~~~~
   + \frac12 \int_{\Omega_0} f'(\fhi\ee) |\nabla \fhi\ee|^2 
   + \frac{1}{4\epsi^2} \iop | \nabla \fhi\ee|^2
   + \frac{1}{4\epsi^2} \iom | \nabla \fhi\ee|^2\\
 \label{energyee4}
  & \mbox{}~~~~~
   \lesssim \io (1 - \fhi\ee)^+ S\ee \mu\ee 
   + c
   +  \io F\ee(\fhi\ee) 
   + c \| \fhi\ee \|_V^2.
\end{align}
It now remains to provide a control of the first summand on the \rhs. 
Using \eqref{S2ee}, we obtain
\begin{equation}\label{co11x}
 \io (1 - \fhi\ee)^+ S\ee \mu\ee 
  = - \io (1 - \fhi\ee^2)^+ \gamma(\cdot,\cdot,\fhi\ee) ( A \fhi \ee + f\ee(\fhi\ee) - \lambda \fhi\ee).
\end{equation}
Now, it is easy to see that, as before,
\begin{equation}\label{co11y}
  - \io (1 - \fhi\ee^2)^+ \gamma(\cdot,\cdot,\fhi\ee) ( A \fhi \ee - \lambda \fhi\ee)
   \le \frac14 \| A\fhi\ee \|^2 + c.
\end{equation}
Hence, it remains to control the quantity
\begin{equation}\label{ee:13}
  - \io (1 - \fhi\ee^2)^+ \gamma(x,t,\fhi\ee) f\ee(\fhi\ee),
\end{equation}
which, however, is readily estimated simply by noting that
$|f\ee(r)| \le |f(r)|$ for $|r| \le 1$ while $(1 - r^2)^+$ is zero for $|r|\ge 1$.
Hence, the argument in \eqref{co11b} basically remains valid up to notational
variations.

\smallskip

\noindent%
Thanks to the above considerations, we finally obtain
\begin{align}\no
  & \ddt \bigg[ \frac12 \| \nabla \fhi\ee \|^2
   + \io \Phi\ee(\fhi\ee)  
   + \frac12 \io F\ee(\fhi\ee) \bigg]
   + \frac12 \| \nabla \mu\ee \|^2
   + \frac12 \| A \fhi\ee \|^2
   + \| \vu\ee \|^2\\
  \no
 & \mbox{}~~~~~~~~~~
   + \frac12 \int_{\Omega_0} f'(\fhi\ee) |\nabla \fhi\ee|^2 
   + \frac{1}{4\epsi^2} \iop | \nabla \fhi\ee|^2
   + \frac{1}{4\epsi^2} \iom | \nabla \fhi\ee|^2\\
 \label{energyee5}
  & \mbox{}~~~~~
   \le c +  \io F\ee(\fhi\ee) 
   + c \| \fhi\ee \|_V^2
   \lesssim c + c \Big[ \frac12 \| \nabla \fhi\ee \|^2
   + \io \Phi\ee(\fhi\ee)  
   + \frac12 \io F\ee(\fhi\ee) \Big],
\end{align}
the last inequality following from the nonnegativity of $\Fhi\ee$ and the 
uniform (in $\epsi$) coercivity of $F\ee$.
Hence the Gr\"onwall lemma can be applied to deduce a number of a-priori estimates
independent of the approximation parameter~$\epsi$. 

Let us note in particular that $F\ee \le F$ everywhere; for this reason, the approximate
energy (i.e., the quantity in square brackets) is controlled at the initial time
by means of~\eqref{hp:init}. Hence, we obtain
\begin{align} \label{st:11}
  & \| \fhi\ee \|_{L^\infty(0,T;V)} + \| \fhi\ee \|_{L^2(0,T;H^2(\Omega))} \le c,\\ 
 \label{st:12}
  & \| F\ee(\fhi\ee) \|_{L^\infty(0,T;L^1(\Omega))} \le c,\\
 \label{st:13}
  & \| \vu\ee \|_{L^2(0,T;H)} \le c,\\
 \label{st:14}
  & \| \nabla \mu\ee \|_{L^2(0,T;H)} \le c,\\
 \label{st:15}
  & \int_0^T \int_{\Omega_+(t)} | \nabla \fhi\ee(t) |^2 \,\dit 
  + \int_0^T \int_{\Omega_-(t)} | \nabla \fhi\ee(t) |^2 \,\dit \le c \epsi^2.
\end{align}
We now provide a control of the mean value of $\fhi\ee$. Namely, we prove 
that, at least for $\epsi>0$ small enough, $(\fhi\ee)\OO$ takes values into the 
reference interval $(-1,1)$. As noted in the introduction, the validity 
of this property is tied to the choice of the ``strongly separating'' potential
\eqref{logpot}.

Using \eqref{st:12} and applying Jensen's inequality, we observe that
there exists a constant $C_T$, possibly depending on $T$ but independent
of $\epsi$, such that, for almost every $t\in(0,T)$, there holds
\begin{equation}\label{co21}
   C_T \ge \io F\ee(\fhi\ee(t,x))\,\frac{\dix}{|\Omega|}
    \ge F\ee \bigg( \io \fhi\ee(t,x)\,\frac{\dix}{|\Omega|} \bigg).
\end{equation}
Then, since $F\ee$ is even and coincides with $F$ on the interval
$[-1+\epsi,1-\epsi]$, we may combine \eqref{co21} with the corresponding property
at the initial time \eqref{jen}. Hence, using in particular 
\eqref{hp:init2}, it is not difficult to check that there exists a number
$\delta\in(0,1/4)$, which depends on the initial datum and on the constant 
$C_T$, but is otherwise independent of $\epsi$, such that
\begin{equation}\label{co22}
   - 1 + \delta \le (\fhi\ee(t))\OO \le - 1 + \delta
    \quext{for every }\,\epsi\in(0,\delta)~~\text{and a.e.~}\,t\in(0,T).
\end{equation}
%
%
%
%
The above property permits us to estimate the mean value of $f\ee(\fhi\ee)$ 
by using a procedure which is rather standard for the Cahn-Hilliard system. Namely, we test \eqref{muee}
by $\fhi\ee - (\fhi\ee)\OO$. Then, using also the Poincar\'e-Wirtinger inequality
it is not difficult to obtain
\begin{align}\no
   & \io f\ee(\fhi\ee)(\fhi\ee - (\fhi\ee)\oo )
   + \| \nabla \fhi\ee \|^2 
  = \io \mu\ee(\fhi\ee - (\fhi\ee)\oo )
   + \lambda \io \fhi\ee(\fhi\ee - (\fhi\ee)\oo )\\
  \no
  & \mbox{}~~~~~
  = \io (\mu\ee-\mu\ee)\oo(\fhi\ee - (\fhi\ee)\oo )
   + \lambda \io (\fhi\ee - (\fhi\ee)\oo )^2\\
 \label{co23}
  & \mbox{}~~~~~
  \le c  \| \nabla \fhi\ee \| \big( \| \nabla\mu\ee \| + \| \nabla \fhi\ee \| \big)
   \le c  \big( 1 + \| \nabla\mu\ee \| \big),
\end{align}
where estimate \eqref{st:11} has been used in order to deduce the last inequality. 

Now, reasoning as in \cite[Prop.~A.1]{MZ}, we deduce
\begin{equation}\label{co24}
  \io f\ee(\fhi\ee)(\fhi\ee - (\fhi\ee)\oo )
   \ge \kappa \| f\ee(\fhi\ee) \|_{L^1(\Omega)} - c,
\end{equation}
where the constants $\kappa>0$ (possibly small) and $c>0$ (possibly large) may
depend on $\delta$ but are independent of $\epsi$.
Hence, replacing \eqref{co24} into \eqref{co23}, squaring, and integrating in time,
we infer 
\begin{equation}\label{st:17}
  \| f\ee(\fhi\ee) \|_{L^2(0,T;L^1(\Omega))}^2
   \le c \big( 1 + \| \nabla \mu\ee \|_{L^2(0,T;H)}^2 \big)
   \le c,
\end{equation}
the last inequality following from \eqref{st:14}.

Hence, integrating \eqref{muee} in space, using \eqref{st:17} and performing
standard manipulations, it is not difficult to deduce
\begin{equation}\label{co25}
  \| (\mu\ee)\oo \|_{L^2(0,T)} \le c,
\end{equation}
which, compared with \eqref{st:14}, yields
\begin{equation}\label{st:18}
  \| \mu\ee \|_{L^2(0,T;V)}
   \le c.
\end{equation}

Let us now observe that, from \eqref{st:18} and \eqref{st:14}, using standard
interpolation and three-dimensional embedding inequalities, there follows
\begin{equation}\label{st:11k}
  \| \mu\ee \nabla\fhi\ee \|_{L^2(0,T;L^{3/2}(\Omega))} \le c,
\end{equation}
whence, from \eqref{st:13}, \eqref{backuee}
and Poincar\'e's inequality, we also get
\begin{equation}\label{st:11j}
  \| p\ee \|_{L^2(0,T;W^{1,3/2}(\Omega))} \le c.
\end{equation}
Now, seeing \eqref{muee} as a family of time independent elliptic
problems with monotone nonlinearities and applying techniques
which are now rather well established (we refer the reader
to \cite{GGW} for details), we infer
\begin{align}\label{st:19}
  & \| \fhi\ee \|_{L^4(0,T;H^2(\Omega))} + \| \fhi\ee \|_{L^2(0,T;W^{2,6}(\Omega))}
   \le c,\\
 \label{st:19b}
   & \| f\ee(\fhi\ee) \|_{L^2(0,T;L^6(\Omega))}
   \le c.
\end{align}


\subsection{Passage to the limit}
\label{subsec:limit}

We now prove that that it is possible take the limit $\epsi\to 0$ in the approximation
scheme detailed above, so to obtain existence of a solution to the original 
system \eqref{phase}-\eqref{S2}.

To this aim, we first observe that, thanks to the procedure carried out in the previous
section, estimates \eqref{st:11}-\eqref{st:15}, \eqref{st:18}
and \eqref{st:11j}-\eqref{st:19b}
hold uniformly with respect to the approximating parameter $\epsi$. 
These estimates, as well as the ones that will follow, imply, by
means of standard weak or weak star compactness results, appropriate convergence
properties, which will be implicitly (i.e., not stressing it in the notation)
intended to hold up to the extraction of subsequences.
First of all, we have
\begin{align} \label{lim:11}
  & \fhi\ee\to\fhi \quext{weakly star in }\,\LIV \cap L^4(0,T;H^2(\Omega)) 
         \cap L^2(0,T;W^{2,6}(\Omega)),\\ 
 \label{lim:12}
  & \vu\ee\to\vu \quext{weakly in }\,\LDH,\\ 
 \label{lim:12b}
  & p\ee\to p \quext{weakly in }\,L^2(0,T;W^{1,3/2}(\Omega)),\\ 
 \label{lim:13}
  & \mu\ee\to\mu \quext{weakly in }\,\LDV.
\end{align}
In order to detail an estimate for $\fhi\eet$, we first replace the 
expression \eqref{backuee} for $\vu\ee$ in the system so to rewrite
\eqref{phaseee} in the more natural form
\begin{equation}\label{phaseee3}
  \fhi\eet + A \mu\ee = ( 1 - \fhi\ee )^+ S\ee - \vu\ee \cdot \nabla \fhi\ee.
\end{equation}
Next, combining the first \eqref{lim:11} with \eqref{lim:12}, we deduce
\begin{equation} \label{lim:21-0}
  \| \vu\ee \cdot \nabla \fhi\ee \|_{L^2(0,T;L^1(\Omega))}\le c.
\end{equation}
Then, using \eqref{lim:13} with the properties of $A$ as a bounded
linear operator from $V$ to $V'$, comparing terms in \eqref{phaseee3} 
we readily deduce
\begin{equation}\label{st:16}
  \| \fhi\eet \|_{L^2(0,T;V' + L^1(\Omega))} \le c.
\end{equation}
Combining the above with \eqref{lim:11} and applying the Aubin-Lions lemma
(see, e.g., \cite{Si}), we then 
obtain
\begin{equation}\label{lim:14}
  \fhi\ee\to \fhi \quext{strongly in }\, C^0([0,T];H^{1-\sigma}(\Omega))
  \cap L^2(0,T;W^{2-\sigma,6}(\Omega)),
\end{equation}
for every $\sigma>0$.

Next, we observe that, thanks to \eqref{st:19b}, there exists a function 
$\xi$ such that 
\begin{equation}\label{lim:15}
  f\ee(\fhi\ee) \to \xi \quext{weakly in }\, L^2(0,T;L^6(\Omega)).
\end{equation}
Hence, {\sl a fortiori}, convergence \eqref{lim:14} holds {\sl strongly}\/ in 
$L^2(0,T;H)$, while \eqref{lim:15} holds {\sl weakly} in $L^2(0,T;H)$.
Thus, in view of the fact that $f\ee$ converges to $f$ in the sense of graphs
(or ``$G$-convergence'', cf.~\cite[Def.~3.58]{At}) in $\RR$, and 
consequently the maximal monotone operators generated by $f\ee$ on $L^2(0,T;H)$
converge, in the sense of graphs in $L^2(0,T;H)$, to the maximal monotone operator 
generated by $f$ on $L^2(0,T;H)$, we can apply a standard monotonicity argument 
in the Hilbert space $L^2(0,T;H)$ (see, e.g., \cite{Br} or \cite{Ba})
to deduce that 
\begin{equation}\label{identif}
  \xi(x,t) = f(\fhi(x,t)) \quext{for a.e.~}\,(x,t)\in \Omega\times(0,T).
\end{equation}
We also notice that, thanks to \eqref{lim:12}, the second \eqref{lim:14} and 
Sobolev's embeddings, there holds at least
\begin{equation} \label{lim:21}
  \vu\ee \cdot \nabla \fhi\ee \to \vu \cdot \nabla \fhi
   \quext{weakly in }\, L^1(0,T;H).
\end{equation}
Moreover, using \eqref{lim:11}, \eqref{lim:13} and \eqref{lim:14} 
(recall also \eqref{st:11k}), we infer
\begin{equation} \label{lim:22}
  \mu\ee \nabla \fhi\ee \to \mu \nabla \fhi
   \quext{weakly in }\, L^2(0,T;L^{3/2}(\Omega)).
\end{equation}
Next, note that, by \eqref{identif}, it follows that $\fhi(x,t)\in(-1,1)$ 
for almost every $(x,t)\in\Omega\times(0,T)$. Combining this with the 
Lipschitz continuity of the positive part function, we then obtain
\begin{align} \label{lim:23}
  & ( 1 - \fhi\ee )^+ \to ( 1 - \fhi )^+ = ( 1 - \fhi) 
   \quext{strongly, say, in }\, L^2(0,T;L^2(\Omega)),\\ 
 \label{lim:24}
  & ( 1 + \fhi\ee )^+ \to ( 1 + \fhi )^+ = ( 1 + \fhi) 
   \quext{strongly, say, in }\, L^2(0,T;L^2(\Omega)).
\end{align}
Analogously, the boundedness of $\gamma$ and its Lipschitz continuity with respect to $\fhi$
(cf.~\eqref{hp:gamma}-\eqref{hp:gamma2}), thanks to the dominated convergence theorem,
imply
\begin{equation} \label{lim:25}
  \gamma(x,t,\fhi\ee) \to \gamma(x,t,\fhi) 
   \quext{strongly in }\, L^r((0,T)\times\Omega)
\end{equation}
for every $r\in[1,\infty)$. This implies, in turn, that 
\begin{equation} \label{lim:25b}
  S\ee = - ( 1 + \fhi\ee)^+ \gamma(x,t,\fhi\ee) \to 
   - ( 1 + \fhi) \gamma(x,t,\fhi) =: S
   \quext{strongly in }\, L^r((0,T)\times\Omega),
\end{equation}
still for every $r\in[1,\infty)$. 

%
%
%
%
%
%

Collecting the above relations, it is readily seen that we can take the 
limit as $\epsi\searrow 0$ of the system given by \eqref{phaseee3} with 
\eqref{muee}-\eqref{S2ee}, and get \eqref{phase}-\eqref{S2} in the limit. 
In particular, the limit of all nonlinear terms is correctly identified.
Finally, we may notice that, in view of \eqref{lim:14},
$\fhi_0 = \fhi\ee|_{t=0} \to \fhi|_{t=0}$ strongly in $H$; hence, the initial
condition \eqref{init} is satisfied. This concludes the proof
of Theorem~\ref{teo:main}.


\subsection{Further remarks}
\label{subsec:rema}

With the existence proof at hand, we would like to give here 
some more comments regarding our assumptions, with particular
reference to the choice of the boundary conditions.

First of all, it is worth observing that the strategy used to derive the 
energy estimate, and in particular the key step (namely, the 
control \eqref{co1b}) would work also in the case of the ``standard'' logarithmic potential
\eqref{logpot0} up to purely technical variations. Indeed, \eqref{co11b}
still holds since $f$ is ``less singular'' in the logarithmic case; 
hence the degenerate character of $(1-\fhi)^2$ balances it with
no need of taking sign conditions on $\gamma$. 
On the other hand, with the logarithmic potential \eqref{logpot0}, the outcome of the 
energy estimate is not sufficient to pass to the limit. Indeed, from the 
bound corresponding to \eqref{st:12} one would no longer be able to deduce 
\eqref{co22} because $F$ takes finite values into the 
{\sl closed} interval $[-1,1]$. When no mass source is present, that would not be
a problem because the value $\fhi\oo$ is conserved;
On the other hand, in our situation, the only other possible
strategy to get a control of $\fhi\oo$ would be that of integrating
\eqref{phase} with respect to space. Nevertheless, as noted in the introduction,
{\sl with the current choice of the boundary conditions}, it seems not possible
to control the last term in \eqref{medie}, whence the mentioned ``inconsistency'' 
phenomenon may occur, leading to ill-posedness
of the system. As noted before, that issue was avoided in \cite{GLRS} by assuming the 
boundary condition \eqref{coupled} and suitably designing the
mass source term. As we consider, instead, the ``strongly singular'' potential 
\eqref{logpot} with the boundary conditions \eqref{bound}, the spatial mean $\fhi\oo$ 
is automatically controlled {\sl from the energy bound}, with no need of 
integrating \eqref{phase} in space.

It is also worth observing that there are other significant choices 
of boundary conditions for which existence of a solution remains,
up to our knowledge, an open issue. To explain this fact,
we go back to the paper \cite{GLRS}, and observe that, there, equation 
\eqref{vu} was in fact replaced by 
\begin{equation}\label{vuglrs}
   \vu = - \nabla q - \fhi \nabla \mu,
\end{equation}
This means that, if $p$ is ``our'' pressure and $q$ is the pressure in \cite{GLRS},
there holds
\begin{equation}\label{pq}
  p = q + \fhi \mu,
\end{equation}
as one can readily see by comparing \eqref{vu} with \eqref{vuglrs}. 
``Incorporating'' the term $\fhi\mu$ into the pressure 
so to obtain a different expression of the Korteweg force 
is a standard procedure, and, indeed, our system \eqref{phase}-\eqref{vu} 
is perfectly equivalent to system~(1.1) of \cite{GLRS}
as far as one looks at the equations on $\Omega$. There is, however, an impact on the
boundary conditions; indeed, in \cite{GLRS} the Dirichlet condition $q=0$ was
considered, while here we are assuming $p=0$, i.e.,
$q + \fhi \mu = 0$. Hence, it is a natural question
to establish whether our arguments may be
extended to the case 
\begin{equation}\label{boundq}
  \dn \mu = \dn \fhi = q = 0 \quext{on }\,\Gamma.
\end{equation}
Note that the above is different both from our \eqref{bound}
and from the choice of \cite{GLRS}, which in our notation reads as
\begin{equation}\label{boundglrs}
  \dn \mu - \fhi \vu\cdot \bn = \dn \fhi = q = 0 \quext{on }\,\Gamma.
\end{equation}
However, it seems that conditions \eqref{boundq} are not straighforward to deal
with, with the main issue arising already at the level of the energy estimate. 
Actually, repeating the computations at the beginning of 
Section~\ref{sec:formal}, i.e.\ combining 
\eqref{energy} with the analogue of \eqref{co12}, 
in the case \eqref{boundq} one would obtain
\begin{equation}\label{energyb}
  \ddt \calE 
   + \| \nabla \mu \|^2 
   + \| \vu \|^2 
   = \io S \mu + \io S q
   - \iga \fhi \mu \vu \cdot \bn,
\end{equation}
where the boundary term has no sign properties and seems very difficult to control
by using trace theorems and embeddings because
only the $H$-norm of $\vu$ appears on the \lhs. Finally, we observe
that the case 
\begin{equation}\label{boundglrs2}
  \dn \mu - \fhi \vu\cdot \bn = \dn \fhi = p = 0 \quext{on }\,\Gamma,
\end{equation}
seems to present mathematical difficulties similar to those occurring 
for \eqref{boundq}; namely, a boundary integral, difficult to be controlled, 
remains on the \rhs\ of the energy inequality.


\section{Additional regularity and uniqueness in dimension two}
\label{sec:rego}

In this part we sketch the proof of Theorem~\ref{teo:rego}. Actually, most of the procedure
will follow the lines of the argument in \cite[Sec.~4]{GLRS}; for this reason we will
only outline the parts where the main differences arise. In fact, all the necessary
variations depend only on the different choice of boundary conditions and to the 
corresponding definition of the pressure (as explained in the previous section), while
at this level the ``strongly separating'' potential plays no role at 
all.

Of course, the main tool in the proof consists in the derivation of additional a-priori estimates,
which, to avoid technicalities, will be carried out in a formal way, i.e.\ by working directly on 
the ``original'' system \eqref{phase}-\eqref{S2} without referring to the approximation.
In particular, we will take advantage of all the information coming from the previous
estimates, including the uniform bound $|\fhi|\le 1$. Actually, to make the procedure completely
rigorous one could proceed along the lines of \cite[Sec.~4]{GLRS} where this issue is discussed in 
detail. Note in particular that the approximation scheme should also be refined,
possibly operating some additional regularization of the initial data.

That said, we first observe that formulas \mbox{\cite[(4.39)-(4.41)]{GLRS}}
remain valid also in the current setting. In our notation, these properties,
which are a direct consequence of the energy estimate, read
\begin{align} \label{glrs:11}
  & \| \mu \|_{V} \le c \big( 1 + \| \nabla \mu \| \big),\\ 
 \label{glrs:12}
  & \| \fhi \|_{H^2(\Omega)}^2 \le c \big( 1 + \| \nabla \mu \| \big),\\ 
 \label{glrs:13}
  & \| \fhi \|_{W^{2,r}(\Omega)} \le c_r \big( 1 + \| \nabla \mu \| \big),
   \quext{for every }\,r\in[1,\infty),
\end{align}
where the constants $c$ (or $c_r$ in the last case, exploding as $r\nearrow\infty$),
only depend on quantities 
that have already been estimated uniformly with respect both to the approximation
parameter and to the time variable. 

Then, the key tool we use in order to obtain higher order 2D estimates is the following
Brezis-Gallouet-Wainger inequality (see, e.g., \cite{Eng}):
\begin{equation}\label{BGW}
  \| f \|_{L^\infty(\Omega)} 
   \le C \| f \|_{V} \log^{1/2} \big ( e + \| f \|_{W^{1,r}(\Omega)} \big)
   + C,
\end{equation}
valid for every $r>2$ and $f\in W^{1,r}(\Omega)$, with $C>0$ also depending on $r$. Of course, 
thanks to Sobolev's embeddings, the $W^{1,r}$-norm may be replaced 
by the $H^2$- one when $f\in H^2(\Omega)$.

That said, we start detailing the regularity estimates. Basically, we
need to test \eqref{phase} by $\mu_t$; 
moreover, we take the time derivative of \eqref{mu} and test it by $\fhi_t$. Then, comparing
the resulting relations, we readily obtain
\begin{equation}\label{re:11}
  \frac12 \ddt \| \nabla \mu \|^2 
   + \| \nabla \fhi_t \|^2 
   + \io f'(\fhi) \fhi_t^2
   - \io (1 - \fhi) S \mu_t
  = \lambda \| \fhi_t \|^2 
   - \io  \mu_t \vu \cdot \nabla\fhi.
\end{equation}
Taking the time derivative of \eqref{vu} and testing it by $\vu$, we also infer
\begin{equation}\label{re:12}
  \frac12 \ddt \| \vu \|^2 
   - \io p_t S 
  = \io \mu_t \vu \cdot \nabla \fhi 
   + \io \mu \vu \cdot \nabla \fhi_t.
\end{equation}
Next, we integrate by parts in time the last terms on the \lhs s of \eqref{re:11}
and \eqref{re:12}. Then, we combine the two obtained relations and notice that
two terms cancel out. Using also the monotonicity of $f$ (which permits us to 
neglect a further positive quantity), we then deduce
\begin{align}\no
  & \ddt \Big[ \frac12 \| \nabla \mu \|^2 
   + \frac12 \| \vu \|^2 
   - \io (1 - \fhi) S \mu
   - \io p S \Big]
   + \| \nabla \fhi_t \|^2 \\
 \label{re:13}
  & \mbox{}~~~~~   
   = \lambda \| \fhi_t \|^2 
   - \io \mu \ddt \big( (1 - \fhi) S \big)
   - \io p S_t
   + \io \mu \vu \cdot \nabla \fhi_t.
\end{align}
This relation basically corresponds to formula (4.56) of \cite{GLRS}. The main
difference is provided by the last term on the \rhs, which in the setting 
of \cite{GLRS} has a different expression, namely
\begin{equation}\label{term:glrs}
  - \io \fhi_t \vu \cdot \nabla \mu,
\end{equation}
due to the different choice of the pressure and of the Korteweg term
taken there.

On the other hand, our situation is in several aspects simpler compared to 
\cite{GLRS}. Actually, here $\mu$ satisfies a no-flux boundary condition, while 
in \cite{GLRS} the boundary conditions for $\mu$ and $\vu$ were coupled
and consequently more complicated to deal with.
Indeed, applying elliptic regularity results to \eqref{phase} 
and recalling \eqref{glrs:11}, it is easy to deduce
\begin{align}\no
  \| \mu \|_{H^2(\Omega)}
   & \le c \big( \| \mu \| + \| \Delta \mu \| \big)
   \le c \big( 1 + \| \nabla \mu \| + \| \fhi_t \| 
    + \| (1 - S) \fhi \| + \| \vu \cdot \nabla \fhi \| \big)\\
 \label{re:14}
  & \le c \big( 1 + \| \nabla \mu \| + \| \fhi_t \| 
    + \| \vu \| \| \nabla \fhi \|_{L^\infty(\Omega)} \big).
\end{align}
This relation, whose analogue in \cite{GLRS} has a much longer proof,
will play a key role in the sequel.
%
%
%
%

To proceed, we would like to provide a control of the \rhs\ of \eqref{re:13}.
To this aim, we need some preparatory computations. First of all, in
order to get a tractable differential inequality, as
in \cite{GLRS} we denote by $\calH$ the quantity under time derivative in 
\eqref{re:13}. Using the uniform boundedness of $S$ and $\fhi$, \eqref{glrs:11}, 
the first \eqref{st:11}, a suitable version of Poincar\'e's inequality,
and Young's inequality, we may then notice that
\begin{align}\no
  \calH & = \Big[ \frac12 \| \nabla \mu \|^2 
   + \frac12 \| \vu \|^2 
   - \io (1 - \fhi) S \mu
   - \io p S \Big] \\
 \no
  & \le \frac12 \| \nabla \mu \|^2 
   + \frac12 \| \vu \|^2 
   + c \| \mu \|
   + c \| p \|_{L^1(\Omega)} \\
 \no
  & \le c \| \nabla \mu \|^2 + c
   + \frac12 \| \vu \|^2 
   + c \| \nabla p \|_{L^1(\Omega)}\\
 \no
  & \le c \| \nabla \mu \|^2 + c
   + \frac12 \| \vu \|^2 
   + c \| \vu \|_{L^1(\Omega)}
   + c \| \mu \| \| \nabla \fhi \|\\
 \label{Hsu}
  & \le c \big( 1 + \| \nabla \mu \|^2 
   + \| \vu \|^2 ).
\end{align}
A similar procedure, which we do not detail for the sake of brevity,
gives also the reverse inequality, i.e.,
\begin{equation}\label{Hgiu}
  \calH \ge \frac14 \big( \| \nabla \mu \|^2 + \| \vu \|^2 \big)
   - C,
\end{equation}
where the ``large'' constant $C>0$ only depends on quantities that have 
already been controlled uniformly in time in the previous estimates.

Let us now observe that, by \eqref{phase} with the Poincar\'e-Wirtinger inequality,
there follows
\begin{align}\no
  \| \fhi_t \| 
   & \le \| \fhi_t - (\fhi_t)\oo \| + c | (\fhi_t)\oo |
     \le c \| \nabla \fhi_t \| + c | (\fhi_t)\oo |\\
 \no
  & \le c \| \nabla \fhi_t \| + c \| (1 - \fhi) S \|_{L^1(\Omega)}
      + c \| \vu \cdot \nabla\fhi \|_{L^1(\Omega)}\\
 \label{fhit:11}
  & \le c \| \nabla \fhi_t \| + c
      + c \| \vu \| \| \nabla\fhi \| 
    \le c \big( \| \nabla \fhi_t \| + \| \vu \| \big).
\end{align}
%
%
%
%
%
%
Then, on account of the above considerations, and using in particular 
\eqref{glrs:11}-\eqref{glrs:13}, \eqref{BGW}, \eqref{re:14}, 
\eqref{Hsu}-\eqref{Hgiu}, and 
\eqref{fhit:11}, the last term in \eqref{re:13} can be estimated as follows:
\begin{align}\no
  \io \mu \vu \cdot \nabla \fhi_t
   & \le \| \mu \|_{L^\infty(\Omega)} \| \vu \| \| \nabla \fhi_t \|\\
 \no   
   & \stackrel{\eqref{BGW},\eqref{glrs:11}}{\le} c \big( 1 + \| \nabla \mu \| \big) \log^{1/2} \big( e + \| \mu \|_{H^2(\Omega)} \big)
    \| \vu \| \| \nabla \fhi_t \|\\
 \no   
   & \stackrel{\eqref{Hgiu}}{\le} c ( 1 + \calH ) \log^{1/2} \big( e + \| \mu \|_{H^2(\Omega)} \big)
    \| \nabla \fhi_t \|\\
 \no   
   & \stackrel{\eqref{re:14}}{\le} c ( 1 + \calH ) \log^{1/2} \big( e + \| \nabla \mu \| + \| \fhi_t \| 
    + \| \vu \| \| \nabla \fhi \|_{L^\infty(\Omega)} \big)
    \| \nabla \fhi_t \|\\
 \no   
   & \stackrel{\eqref{fhit:11},\eqref{Hgiu}}{\le} c ( 1 + \calH ) \log^{1/2} \big( e + \calH^{1/2} 
    + \| \nabla \fhi_t \| 
    + \calH^{1/2} \| \fhi \|_{W^{2,3}(\Omega)} \big)
    \| \nabla \fhi_t \|\\
 \no   
   & \stackrel{\eqref{glrs:13}\eqref{Hgiu}}{\le} c ( 1 + \calH ) \log^{1/2} 
   \big( e + \| \nabla \fhi_t \| + \calH \big)
    \| \nabla \fhi_t \|\\ \no   
 \no   
   & \le c ( 1 + \calH ) \log^{1/2} ( e + \calH )  \| \nabla \fhi_t \| 
    + c ( 1 + \calH ) \log^{1/2} \big( e + \| \nabla \fhi_t \| \big) 
  \| \nabla \fhi_t \|\\
 \label{re:17} 
   & =: J_1 + J_2.
\end{align}
Let us now go back to \eqref{re:13} and observe first that, by interpolation,
\begin{align}\no
  \| \fhi_t \| 
   & \le \| \fhi_t - (\fhi_t)\oo \| + | (\fhi_t)\oo |
   \le \| \fhi_t - (\fhi_t)\oo \|_{V'}^{1/2} \| \fhi_t - (\fhi_t)\oo \|_{V}^{1/2} 
    + | (\fhi_t)\oo |\\
 \label{co:powi} 
   & \le c \| \fhi_t \|_{V'}^{1/2} \| \nabla \fhi_t \|^{1/2} 
    + c \| \fhi_t \|_{V'}
   \le \epsilon \| \nabla \fhi_t \| + c_\epsilon \| \fhi_t \|_{V'},
\end{align}
for every ``small'' $\epsilon>0$ and correspondingly ``large''
$c_\epsilon>0$.

Then, using the above relation (with $\epsilon$ taken small enough,
which generates the constant $7/8$ below)
it is very easy to provide a control of 
the first three terms on the \rhs\ of \eqref{re:13},
all of which basically depend on the $L^2$-norm of $\fhi_t$
(see \cite{GLRS} for details). Managing the last summand in \eqref{re:13}
by means of \eqref{re:17}, it is then not difficult to obtain
\begin{equation}\label{re:13b}
  \ddt \calH
   + \frac78 \| \nabla \fhi_t \|^2
   \le c ( 1 + \calH )^2 + J_1 + J_2.
\end{equation}
Moreover, using Young's inequality it is apparent that
\begin{equation}\label{re:18}
  J_1 = c ( 1 + \calH ) \log^{1/2} ( e + \calH ) \| \nabla \fhi_t \| 
   \le \frac18 \| \nabla \fhi_t \|^2 + c ( 1 + \calH )^2 \log ( e + \calH ).
\end{equation}
The control of $J_2$ is a bit more delicate. To achieve it, we 
consider the convex function 
\begin{equation}\label{psi:11}
  \psi : [0,\infty) \to \RR, \qquad 
   \psi(r) = r^2 \log (e + r)
\end{equation}
and we would like to provide some estimate near infinity for its convex conjugate  
\begin{equation}\label{psi:12}
  \psi^*(s) = \sup_{r\ge0} \{ rs - \psi(r) \}.
\end{equation}
Actually, a simple computation permits us to check that,
when $s$ is assigned large enough, the supremum is assumed
at $r=r_0$ such that
\begin{equation}\label{psi:13}
  r_0=(\psi')^{-1}(s), \quext{i.e., }\,
   s = \psi'(r_0) = 2r_0 \log (e + r_0) + \frac{r_0^2}{e + r_0}
   \sim 2r_0 \log (e + r_0),
\end{equation}
where the last approximate equality holds as $s$, hence $r_0$, is large enough. In particular, 
using the relation $\log ( e + r_0) \le r_0/2$, also holding for $r_0$ large
enough, we may estimate
\begin{equation}\label{psi:14}
  \log ( e + s ) \lesssim \log \big( e + 2r_0 \log (e + r_0) \big)
   \le \log \big( e + r_0^2 \big) \le \log ( e + r_0 )^2
   = 2 \log (e + r_0).
\end{equation}
Hence, replacing the above value of $s$ in
\eqref{psi:12}, we deduce
\begin{align}\no
  \psi^*(s) & \sim 2 r_0^2 \log ( e + r_0 ) - \psi (r_0)
   = r_0^2 \log ( e + r_0 )\\  
 \label{psi:15}
  & \sim \frac{s^2}{4 \log^2(e+r_0) }\log ( e + r_0 )
    = \frac{s^2}{4 \log(e+r_0) } 
    \lesssim \frac{s^2}{2 \log(e+s) }.
\end{align}
Then, we use the Fenchel inequality
\begin{equation}\label{fench}
  sr \le \psi(r) + \psi^*(s)
\end{equation}
with the choices
\begin{equation}\label{fench2}
  r = c ( 1 + \calH ), \qquad 
   s = \| \nabla \fhi_t\| \log^{1/2} \big( e + \| \nabla \fhi_t \| \big),
\end{equation}
where the constant $c$ is the same as in the definition
of $J_2$ in~\eqref{re:17}. It is then easy to verify that
\begin{equation}\label{fench3}
  \psi(r) \le c (1 + \calH)^2 \log (e + \calH). 
\end{equation}
On the other hand, by \eqref{psi:15}, we deduce 
\begin{align}\no
  \psi^*(s) & \lesssim \frac{s^2}{2 \log(e+s) }
   \sim \frac{ \| \nabla \fhi_t \|^2 \log \big( e + \| \nabla \fhi_t \| \big) }{%
   2 \log\big( e + \| \nabla \fhi_t\| \log^{1/2} ( e + \| \nabla \fhi_t \| ) \big) }\\
 \label{re:19}
   & \lesssim \frac{ \| \nabla \fhi_t \|^2 \log \big( e + \| \nabla \fhi_t \| \big) }{%
   2 \log\big( e + \| \nabla \fhi_t\| \big) }
    \lesssim \frac12 \| \nabla \fhi_t \|^2.
\end{align}
Consequently, by \eqref{re:18}, \eqref{fench3} and \eqref{re:19},
\eqref{re:13b} reduces to
\begin{equation}\label{re:13x}
  \ddt ( e + \calH )
   + \frac14\| \nabla \fhi_t \|^2
   \lesssim c (e + \calH)^2 \log (e + \calH). 
\end{equation}
Hence, in view of the fact that, 
as a consequence of the energy estimate, $\calH\in L^1(0,T)$
(cf.\ in particular \eqref{st:13} and \eqref{st:18}), we may
apply a generalized version of Gr\"onwall's lemma to the above relation.
This yields the regularity properties
\begin{align} \label{rg:11}
  & \nabla\mu \in \LIH , \qquad \vu \in \LIH,\\
 \label{rg:12}
  & \fhi_t \in \LDV,
\end{align}
provided that the value of $\calH$ is finite at the initial time $t=0$, which, 
at least formally, corresponds to the conditions
\begin{equation}\label{umu0}
   \mu_0:=\mu|_{t=0} \in V, \qquad \vu_0:=\vu|_{t=0} \in H.
\end{equation}
Actually, due to the quasi-stationary character of the model, the 
above properties have to be deduced from the regularity assumed
on $\fhi_0$, which is the sole initial datum associated with the 
system. Hence, we now verify that hypotheses \eqref{init:reg1}-\eqref{init:reg2}
imply \eqref{umu0}. Of course, this is just a formal check: indeed,
to make the procedure fully rigorous one should
properly intervene on the approximation argument
at the price of further technicalities.
That said, we notice that assumption \eqref{init:reg2},
with the position \eqref{defi:mu0}, corresponds exactly
to the first \eqref{umu0}. Hence we just need to show that 
$\vu_0 \in H$. To see this, we consider the elliptic
problem (formally) associated to \eqref{vu} at the initial time,
namely
\begin{equation}\label{def:p0}
  B p_0 = S|_{t=0} - \dive (\mu_0 \nabla \fhi_0 ),
   \qquext{so that }\, \vu_0 = - \nabla p_0 + \mu_0 \nabla \fhi_0.
\end{equation}
Then, we notice that the initial value of $S$ depends on $\fhi_0$ only, and it 
is easy to check that it lies in $H$ as a consequence of \eqref{init:reg1} and
\eqref{hp:gamma}. Moreover, we have 
\begin{equation}\label{p0:11}
  \| \mu_0 \nabla \fhi_0 \| 
   \le \| \mu_0 \|_{L^4(\Omega)} \| \nabla \fhi_0 \|_{L^4(\Omega)}
   \le c \| \mu_0 \|_{V} \| \fhi_0 \|_{H^2(\Omega)}
   < +\infty,
\end{equation}
the last inequality following from \eqref{init:reg1}-\eqref{init:reg2}.
Hence, $\dive (\mu_0 \nabla \fhi_0 ) \in H^{-1}(\Omega)$ and,
by elliptic regularity results, $p_0 \in H^1_0(\Omega)$, which
in turn implies $\nabla p_0 \in H$ and, by \eqref{p0:11} again,
$\vu_0\in H$, as desired. Hence the initial value of $\calH$ is 
finite, which implies \eqref{rg:11}-\eqref{rg:12}.

To conclude the proof, we need to show the regularity properties  
\eqref{rego:fhis}-\eqref{rego:vus}. Actually, a part of them follows from
\eqref{rg:11}-\eqref{rg:12}. To prove the missing ones, we proceed
along the lines of \cite{GLRS} with some small variations deriving from
the different expression of the Korteweg term. 

First of all, we combine again \eqref{co23} and \eqref{co24}. Then, by
means of the improved time regularity of $\nabla \mu$ in \eqref{rg:11},
it is not difficult to improve \eqref{st:17} and \eqref{st:18} 
as follows:
\begin{equation}\label{rg:13}
  f(\fhi) \in L^\infty(0,T;L^1(\Omega)), \qquad \mu \in \LIV.
\end{equation}
With this at hand, noting that in 2D one has 
$V\subset L^r(\Omega)$ for every $r\in[1,\infty)$, 
and applying elliptic regularity results to the analogue
of \eqref{muee}, we deduce \eqref{rego:ffhis} and the last
\eqref{rego:fhis}.

Next, we observe that $p$ solves the time dependent
family of elliptic problems
\begin{equation}\label{ell:p}
  B p = S - \nabla \mu \cdot \nabla \fhi
   - \mu \Delta \fhi,
\end{equation}
where it is easy to check that the \rhs\ lies in $L^\infty(0,T;H)$ thanks
to \eqref{rg:13} and the last \eqref{rego:fhis} (which implies, in 
particular, that $\nabla\fhi$ is bounded in the uniform norm).
From standard elliptic regularity
results we then deduce \eqref{rego:ps}. Now, reasoning on \eqref{phaseteo}
and using the boundedness of $A$ as a linear operator from $V$ to $V'$, 
it is not difficult to get the first \eqref{rego:fhis}. From \eqref{rego:fhis}
and interpolation, we then also obtain $\fhi_t\in L^4(0,T;H)$. In turn, seeing 
\eqref{phaseteo} as an elliptic problem for $\mu$, this yields the
second \eqref{rego:mus}. Finally, a direct check shows that
$\mu\nabla\fhi\in \LIV$. Combining this with \eqref{rego:ps}, we then
deduce \eqref{rego:vus}, which concludes the proof of the regularity
part.

\smallskip

We now sketch the proof of uniqueness, which, as noted above, may be 
carried out basically following the lines of \cite{GLRS}. To this aim, we consider two strong 
solutions $(\fhi_1,\vu_1,\mu_1,p_1)$ and $(\fhi_2,\vu_2,\mu_2,p_2)$ originating from the 
same initial condition $\fhi_0$. We then set
\begin{equation}\label{defi:diff}
  \fhi:=\fhi_1-\fhi_2, \quad \vu:=\vu_1-\vu_2, 
   \quad p:=p_1-p_2, \quad \mu:=\mu_1-\mu_2,
\end{equation}
so that the quadruplet $(\fhi,\vu,p,\mu)$ turns out to solve the
system
\begin{align}\label{phasediff}
  & \fhi_t + \vu_1 \cdot \nabla \fhi + \vu \cdot \nabla \fhi_2
   = \Delta \mu + ( 1 - \fhi_1 ) S_1 - ( 1 - \fhi_2 ) S_2, \\
 \label{mudiff}
  & \mu = - \Delta \fhi + f(\fhi_1) - f(\fhi_2) - \lambda \fhi,\\
 \label{Sdiff}
  & \dive \vu = S_1 - S_2, \\
 \label{vudiff}
  & \vu = - \nabla p + \mu_1 \nabla \fhi + \mu \nabla \fhi_2,
\end{align}
where we have also set $S_i=-(1+\fhi_i)\gamma(x,t,\fhi_i)$, for $i=1,2$.

We then test \eqref{vudiff} by $\vu$, \eqref{phasediff} by $\mu$, 
and \eqref{mudiff} by $\fhi_t$.
Using standard tools, it is then not difficult to arrive at 
\begin{align}\no
  & \ddt\Big[ \frac12 \| \nabla \fhi \|^2 + \frac12 \io \ell(\fhi_1,\fhi_2) | \fhi |^2 \Big]
   + \| \vu \|^2
   + \| \nabla \mu \|^2\\
 \no
  & \mbox{}~~~~~
   = \lambda \io \fhi \fhi_t 
   - \io \mu \vu_1 \cdot \nabla \fhi
   + \io \mu_1 \vu \cdot \nabla \fhi
   + \frac12 \io |\fhi|^2 \ddt \ell(\fhi_1,\fhi_2)\\
 \label{uniq:11}
  & \mbox{}~~~~~~~~~~
   + \io p (S_1 - S_2)
   + \io \big( ( 1 - \fhi_1 ) S_1 - ( 1 - \fhi_2 ) S_2 \big) \mu,
\end{align}
where we have set
\begin{equation}\label{defi:ell}
  \ell(\fhi_1,\fhi_2):= \int_0^1
   f'(s \fhi_1 + (1-s) \fhi_2)\,\dis.
\end{equation}
Moreover, we also need to consider the elliptic problem associated with 
\eqref{vudiff}, namely
\begin{equation}\label{pdiff}
  B p = S_1 - S_2 
   - \dive ( \mu_1 \nabla \fhi + \mu \nabla \fhi_2 ).
\end{equation}
Testing \eqref{pdiff} by $B^{-1} p$, and using standard embedding inequalities, 
it is then easy to deduce
\begin{align}\no
  \| p \|^2 
   & = \io ( S_1 - S_2 ) B^{-1} p
   + \io (\mu_1 \nabla \fhi + \mu \nabla \fhi_2 )\cdot \nabla B^{-1} p\\
 \no
   & \le \frac14 \| p \|^2 + c \| \fhi \|^2
   + \| \mu_1 \nabla \fhi + \mu \nabla \fhi_2 \| \| \nabla B^{-1} p \|\\
 \label{pdiff2}
   & \le \frac12 \| p \|^2
   + c \big( 1 + \| \mu_1 \|_{H^2(\Omega)}^2 \big) \| \fhi \|^2_V
   + c \| \mu \|^2,
\end{align}
where we also used that $\nabla\fhi_2$ is bounded in the uniform norm as a consequence of 
the last \eqref{rego:fhis}.

Then, we test \eqref{mudiff} by $\mu$ to control the last term (note that this will
be used also to treat the other terms depending on $\mu$
on the \rhs\ of \eqref{uniq:11}).
Standard embedding inequalities permit us to arrive at
\begin{equation}\label{mudiff2}
  \| \mu \|^2 
   \le \epsilon \| \nabla \mu \|^2
    + c_\epsilon \big(1 + \| \ell(\fhi_1,\fhi_2) \|_{L^4(\Omega)}^2 \big) \| \fhi \|^2_V,
\end{equation}
for small $\epsilon>0$ and correspondingly large $c_\epsilon>0$ to be chosen later on.

Now, noticing that $\ell\ge 2$ almost everywhere thanks to \eqref{logpot2},
it is easy to verify that 
\begin{equation}\label{Kdiff}
  \calK := \Big[ \frac12 \| \nabla \fhi \|^2 + \frac12 \io \ell(\fhi_1,\fhi_2) | \fhi |^2 \Big]  
   \ge \frac12 \| \fhi \|_V^2.
\end{equation}
Comparing terms in \eqref{phasediff}, we then see that
the first term on the \rhs\ of \eqref{uniq:11} can be treated as follows:
\begin{align}\no
 \lambda \io \fhi \fhi_t 
  & \le \epsilon \| \fhi_t \|_{V'}^2 + c_\epsilon \| \fhi \|_V^2\\
 \no
  & \le \epsilon \big\| \vu_1 \cdot \nabla \fhi + \vu \cdot \nabla \fhi_2
   - \Delta \mu - ( 1 - \fhi_1 ) S_1 + ( 1 - \fhi_2 ) S_2\big \|_{V'}^2
      + c_\epsilon \| \fhi \|_V^2\\
 \no
  & \le c \epsilon \| \vu_1 \cdot \nabla \fhi \|^2_{L^{3/2}(\Omega)}
   + c \epsilon \| \Delta \mu \|_{V'}^2
   + c \epsilon \big\| \vu \cdot \nabla \fhi_2 - ( 1 - \fhi_1 ) S_1 + ( 1 - \fhi_2 ) S_2\big \|^2
   + c_\epsilon \| \fhi \|_V^2\\ 
 \no
  & \le c \epsilon \| \vu_1 \|_{L^{6}(\Omega)}^2 \| \nabla \fhi \|^2
   + c \epsilon \| \nabla \mu \|^2
   + c \epsilon \| \vu \|^2 \| \nabla \fhi_2 \|_{L^{\infty}(\Omega)}^2 + c_\epsilon \| \fhi \|_V^2\\
 \label{uniq:12} 
  & \le \frac18 \big( \| \vu \|^2 + \| \nabla \mu \|^2 \big)
    + c \big( 1 + \| \vu_1 \|_{V}^2 \big) \| \fhi \|_V^2, 
\end{align}
where in the last passage we have taken $\epsilon>0$ small enough depending on
the value of the $L^\infty$-norm of $\nabla\fhi_2$, which is a known quantity.

Then, it is readily seen that the remaining terms on the \rhs\ of \eqref{uniq:11}
can be treated without further difficulties (the details are very similar to 
\cite[Sec.~5]{GLRS}),
provided that one can estimate the quantity $\ell$ defined in \eqref{defi:ell}.
This is in a sense the most delicate part because it involves the singular
function $f$ and its derivatives.
Namely, going back to \eqref{uniq:11} and \eqref{mudiff2},
we need to provide a bound for the sum
\begin{equation}\label{stimell}
  \frac12 \io |\fhi|^2 \ddt \ell(\fhi_1,\fhi_2)
  + c_\epsilon \| \ell(\fhi_1,\fhi_2) \|_{L^4(\Omega)}^2 \| \fhi \|^2_V.
\end{equation}
On the other hand, with the strongly separating potential \eqref{logpot}
this task is in fact simpler compared to the standard logarithmic potential
\eqref{logpot0} considered in \cite{GLRS}. 
Indeed, using \eqref{rego:ffhis} and performing direct computations,
one can easily verify that 
\begin{equation}\label{regofffff}
   f^{(j)}(\fhi_i) \in L^\infty(0,T;L^r(\Omega)) \quext{for all }\,
    j\ge 0,~~i=1,2,~~\text{and }\,r\in[1,\infty).
\end{equation}
Then, we just sketch the estimation of the first summand in \eqref{stimell},
the latter one being in fact simpler. Actually, by \eqref{defi:ell},
\begin{equation}\label{ell:t}
  \ddt \ell(\fhi_1,\fhi_2):= \int_0^1
   f''(s \fhi_1 + (1-s) \fhi_2) (s \fhi_{1,t} + (1-s) \fhi_{2,t})\,\dis,
\end{equation}
whence, by standard embeddings,
\begin{equation}\label{stimell2}
  \frac12 \io |\fhi|^2 \ddt \ell(\fhi_1,\fhi_2)
   \le c \big( 1 + \| f''(\fhi_1) \|_{L^4(\Omega)} + \| f''(\fhi_2) \|_{L^4(\Omega)} \big)
     \big( \| \fhi_{1,t} \|_{V} + \| \fhi_{2,t} \|_{V} \big) \| \fhi \|_V^2.
\end{equation}
Then, going back to \eqref{uniq:11} and noting that the remaining terms on the \rhs\
can be controlled in a simple way (see \cite{GLRS} for details),
it is readily seen that the Gr\"onwall lemma can be applied to 
the functional $\calK$ defined by \eqref{Kdiff}, so to obtain uniqueness
of $\fhi$ and (using the latter summand on the \lhs\ of \eqref{uniq:11} together
with \eqref{mudiff2}) of $\mu$. Then, the uniqueness of $p$ is obtained from the 
elliptic problem \eqref{pdiff} and, finally, that of $\vu$ from \eqref{vudiff}.
This concludes the proof of Theorem~\ref{teo:rego}.
\beos\label{dip:con}
 Refining a bit the procedure (and in particular specifying the control of 
 the remaining quantities on the \rhs\ of \eqref{uniq:11}), one may also obtain
 a continuous dependence estimate, so yielding well-posedness of the model in
 the class of ``strong'' solutions in the 2D setting.
\eddos


\section*{Acknowledgments}

\noindent
This research was supported by the Italian Ministry of Education, University and Research
(MIUR): Dipartimenti di Eccellenza Program (2018-2022), Department of Mathematics ``F.~Casorati'', 
University of Pavia. 
The present paper also benefits from the support of the MIUR-PRIN Grant 2015PA5MP7 ``Calculus of Variations''
and of the GNAMPA (Gruppo Nazionale per l'Analisi Matematica, la Probabilit\`a e le loro Applicazioni) 
of INdAM (Istituto Nazionale di Alta Matematica).


\end{document}